\documentclass[12pt,onecolumn,journal]{IEEEtran}
\usepackage[english]{babel}
\usepackage[latin1]{inputenc}
\usepackage{amssymb}
\usepackage[algo2e,noline,boxruled]{algorithm2e}
\usepackage{cite}
\usepackage{amsmath}
\usepackage{graphicx}
\usepackage[caption=false,font=footnotesize]{subfig}
\usepackage{url}
\usepackage{psfrag}
\usepackage[usenames]{color}
\usepackage{balance}
\usepackage{pstricks}
\usepackage{times}
\usepackage{ifpdf}
\usepackage[]{units}

\newcommand{\beq}{\begin{equation}}
\newcommand{\eeq}{\end{equation}}

\newcommand{\baq}{\begin{eqnarray}}
\newcommand{\eaq}{\end{eqnarray}}

\newcommand{\baqm}{\begin{eqnarray*}}
\newcommand{\eaqm}{\end{eqnarray*}}

\newcommand{\barr}{\begin{array}}
\newcommand{\earr}{\end{array}}

\newcommand{\bi}{\begin{itemize}}
\newcommand{\ei}{\end{itemize}}

\newcommand{\bone}[1]{\mathbf 1\{{#1}\}}

\newcommand{\MSE}{\mbox{MSE}}

\newcommand{\transpose}{^{\ensuremath{\mathsf{T}}}}
\newcommand{\given}{\,\vert\,}

\newcommand{\btheta}{{\boldsymbol \theta}}
\newcommand{\bphi}{{\boldsymbol \phi}}

\usepackage{array}

\newcommand{\cS}{\mathcal{S}}

\newcommand{\bbS}{\mathbb{S}}
\DeclareMathOperator*{\argmin}{arg\,min}

\newcommand{\bn}{\binom}

\newenvironment{pf}{{\bf Proof. }}{\hfill $\square$\medskip}
\newtheorem{theorem}{Theorem}[section]
\newtheorem{lemma}[theorem]{Lemma}

\newcommand{\techreport}[2]{#1}
\newcommand{\fixme}[1]{#1}

\title{On Set Size Distribution Estimation and the Characterization of Large Networks via Sampling}
\author{Fabricio Murai, Bruno Ribeiro, Don Towsley, and Pinghui Wang}

\author{
\IEEEauthorblockN{Fabricio Murai\IEEEauthorrefmark{1}, Bruno Ribeiro\IEEEauthorrefmark{1}, Don Towsley\IEEEauthorrefmark{1}, and Pinghui Wang\IEEEauthorrefmark{2}} \\
\small{\IEEEauthorblockA{\IEEEauthorrefmark{1}Computer Science Department\\
University of Massachusetts Amherst\\ Amherst, MA 01003\\
Email: \{fabricio,ribeiro,towsley\}@cs.umass.edu} \\
\IEEEauthorblockA{\IEEEauthorrefmark{2} State Key Lab for Manufacturing Systems\\
Xi'an Jiaotong University\\
Xi'an P.R.China\\
Email: phwang@sei.xjtu.edu.cn}\\
\vspace{1.0cm}\large{Technical Report UM-CS-2012-023v2}}
}

\begin{document}
\maketitle

\begin{abstract}
In this work we study the set size distribution estimation problem,
where elements are randomly sampled from a collection of non-overlapping sets
and we seek to recover the original set size distribution from the samples. 
This problem has applications to capacity planning, network theory, among other areas.
Examples of real-world applications include characterizing in-degree distributions in large graphs and uncovering TCP/IP flow size
distributions on the Internet.
We demonstrate that it is hard to estimate the original set size distribution.
The recoverability of original set size distributions presents a sharp threshold
with respect to the fraction of elements that remain in the sets. If this
fraction remains below a threshold, typically half of the elements in power-law
and heavier-than-exponential-tailed distributions, then the original set size
distribution is unrecoverable.
We also discuss practical implications of our findings.
\end{abstract}


%
%

\section{Introduction}

Networks are increasingly large and complex, posing tremendous challenges to
their characterization in the wild. Characterizing network structure (e.g.\
        degree distribution), network traffic flows
      (e.g.\ TCP/IP flow sizes in communication networks), node labels
      (e.g.\ group memberships), is usually impossible without resorting to sampling
      due to the size and scale of current networks.
Practitioners often sample networks to estimate their characteristics.
Many problems in network characterization through sampling can be mapped into
the class of set size distribution estimation problems.
The set size distribution estimation problem is stated as follows. Consider a collection of
non-overlapping sets whose elements are probabilistically sampled. 
The problem is to estimate the original (pre-sampling) set size distribution
based on the samples.

Set size distribution estimation has several applications.
One example of particular interest
is the estimation of in-degree distributions of on-line social networks,
where nodes represent people and a directed edge represents, for instance,
one or more messages exchanged between two pairs of nodes.
By monitoring message exchanges one samples a fraction of the edges.
Using these samples we want to estimate the in-degree or out-degree distribution of nodes.
The set size distribution problem also manifests itself in other areas, including Internet traffic monitoring, 
e.g., estimating the size distribution (in packets) of TCP/UDP flows~\cite{Duffield}, and in next generation Internet capacity planing, 
such as estimating the number of copies of a movie in a CDN of next-generation routers.
Fortunately, simple maximum likelihood~\cite{Duffield} or Bayesian-style estimators exist, even when we are unable to observe sets without observed elements.

Despite the importance of characterizing set size distributions, to the best of our knowledge no deep
analysis of set size distribution estimation exists in the literature.
We fill this gap and {\em show that set size distribution estimation exhibits intriguing abnormal statistical properties.}
To best illustrate our results, consider the estimation of in-degree distributions of arbitrarily large power-law graphs.
We prove that if less than 50\% of the edges are observed then the output of {\em any estimator} (be it frequentist or Bayesian) will be as truthful to the original in-degree distribution as a set of random numbers between zero and one.
Moreover, when nodes without sampled incoming edges are unobservable, even a first order metric like average degree is subject to the same threshold behavior, i.e., sampling less than 50\% of all incoming edges impedes the estimation of in-degree averages.
The latter result seemly defies intuition. 
We prove these and other results in the general setting of sets with arbitrary set size distributions.
In what follows we give an overview of our contributions.

\subsection{General Observations}
In this work we uncover intriguing  set size distribution estimation properties, including:
\begin{itemize}
\item {\em A (finite) increase in
samples may result in {\em no} reduction in estimation errors.}
\end{itemize}
Unlike estimation problems such as election
polls, where a sufficient increase in samples always
results in increased accuracy, we show, paradoxically, that in the set size distribution estimation problem
an increase in samples may, in practice, result in no increase in accuracy.
Section~\ref{sec:results} unveils the root cause of this odd behavior and
explains when it can be avoided. Another interesting property is:
\begin{itemize}
\item {\em In networks with large set sizes (e.g., nodes with large degrees) and power-law set size distributions (in fact our results hold for any heavier-than-exponential distributions), randomly sampling less than 50\% of set elements (e.g., edges of a node)
provides {\em almost no} information about the set size distribution or the average set size. 
However, in networks with sub-exponential set size distributions, accurate set
size distributions estimation is always possible.}
\end{itemize}
The above observation is interesting because power-laws
have more tail probability mass and, thus, large sets are more likely to have sampled elements than in sub-exponential tails.
However, and despite this, we show that if less than $50\%$ of elements are sampled, then
estimates of power-laws distributions (more precisely, any heavier-than-exponential distribution) are significantly less accurate than the estimates obtained from sub-exponential distributions.
Our work also provides a host of equally puzzling observations, fully and formally presented in Section~\ref{sec:results}.

%
%
%

\subsection{Outline}
Our paper is organized as follows.
In Section~\ref{sec:motivation} we conduct experiments on
the indegree distribution estimation with real data.
Section~\ref{sec:model} presents the sampling and estimation models.
Section~\ref{sec:results} presents our theoretic results.
Section~\ref{sec:discussion} presents our discussion section where we analyze problems that field analysts are likely to face in practice, highlighting common mistakes made in the literature and how to avoid them.
Finally Section~\ref{sec:conclusion} presents the conclusions and related work.

%
%

\section{Estimation with Real Data}\label{sec:motivation}

In this section, we experiment with one particular application of the set size distribution
problem: the estimation of the in-degree distribution of a network.
Consider the Enron dataset, that describes a network
composed by a group of people
who exchanged emails during a certain period of time.
Here each node represents a person and two people have a directed edge if one has emailed the other. 
The maximum in-degree in this network is 1383.

\begin{figure}[!t]
\centerline{
\subfloat[]{\includegraphics[width
=2.1in]{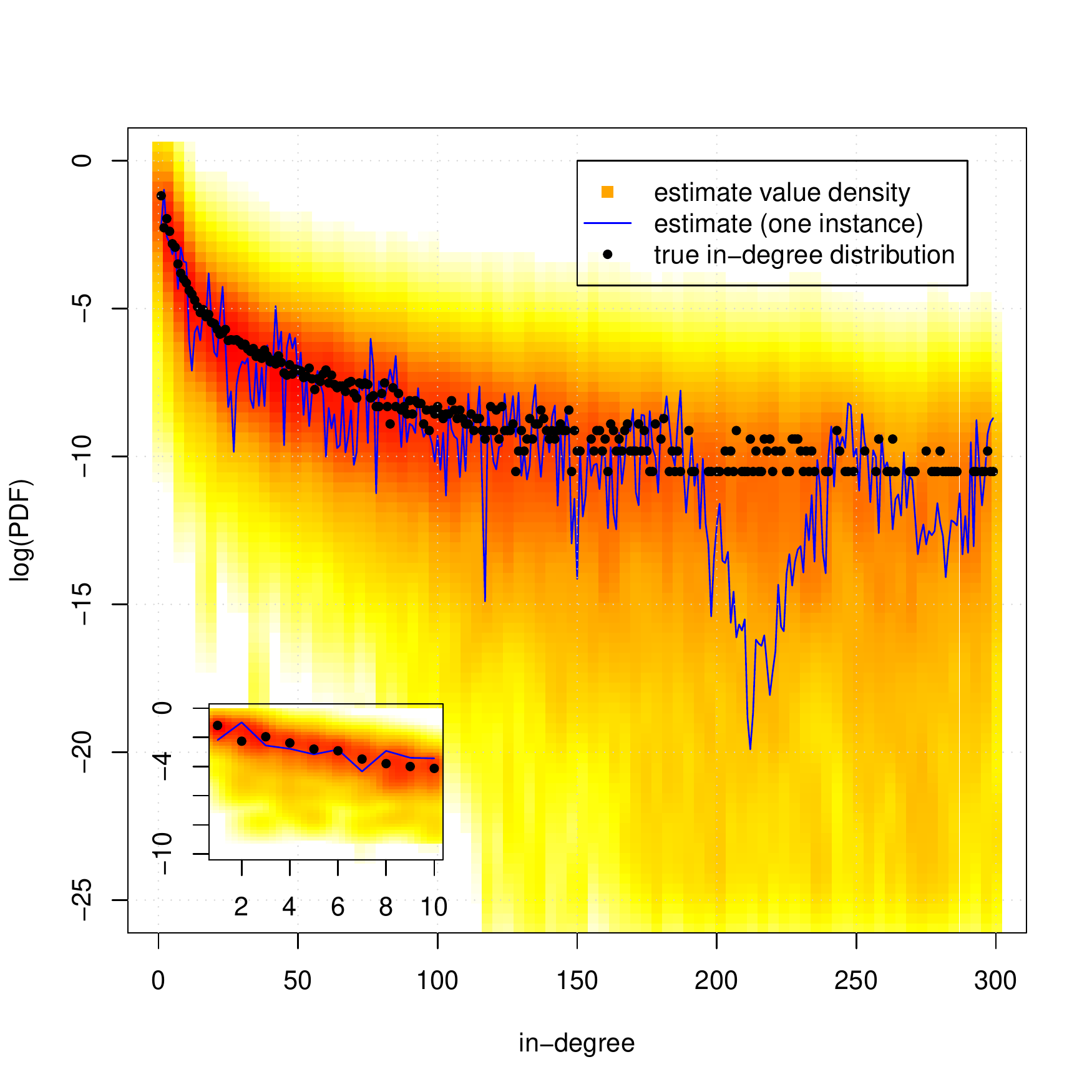}
\label{fig:heatmap_025_10k}}
\subfloat[]{\includegraphics[width
=2.1in]{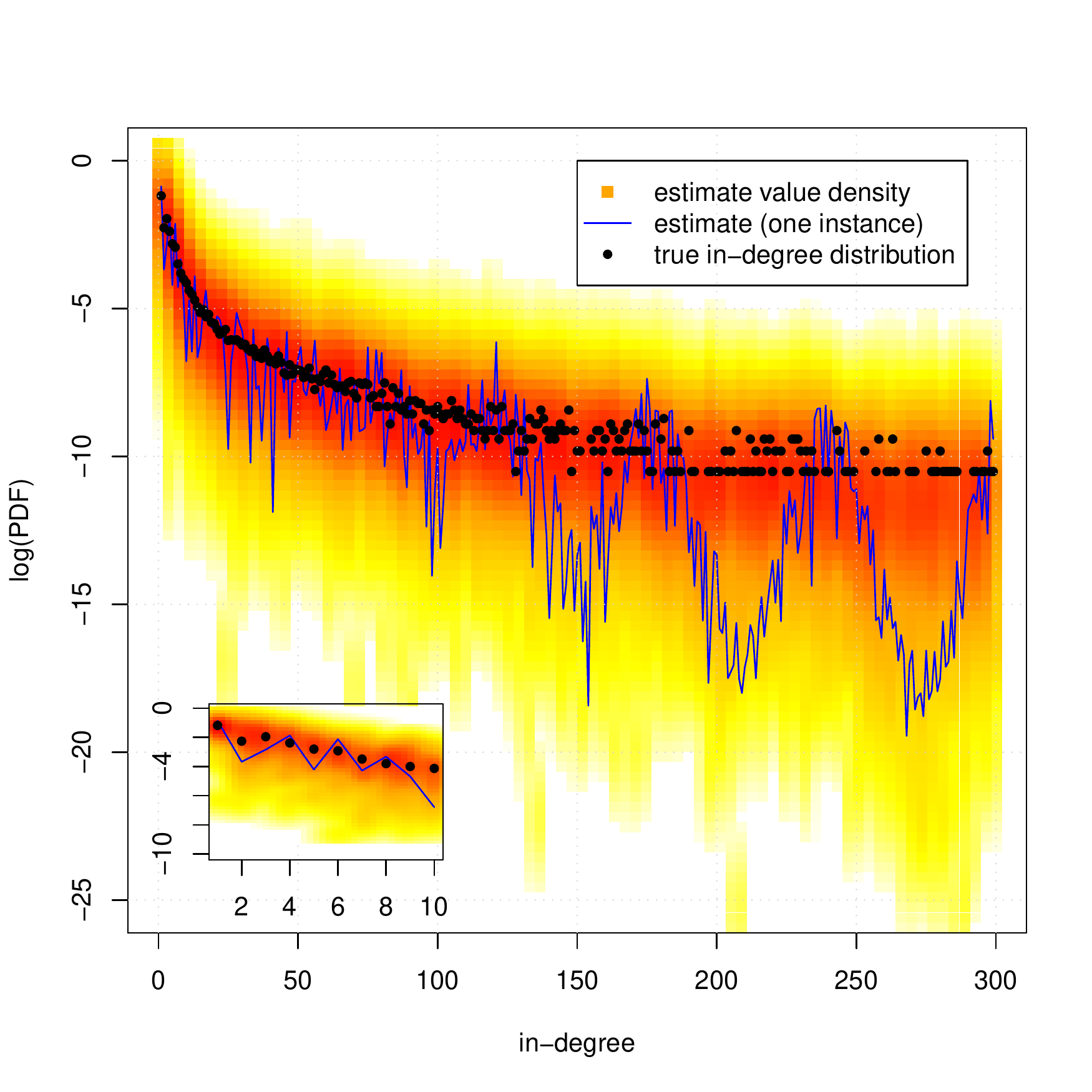}
\label{fig:heatmap_025_50k}}
\subfloat[]{\includegraphics[width
=2.1in]{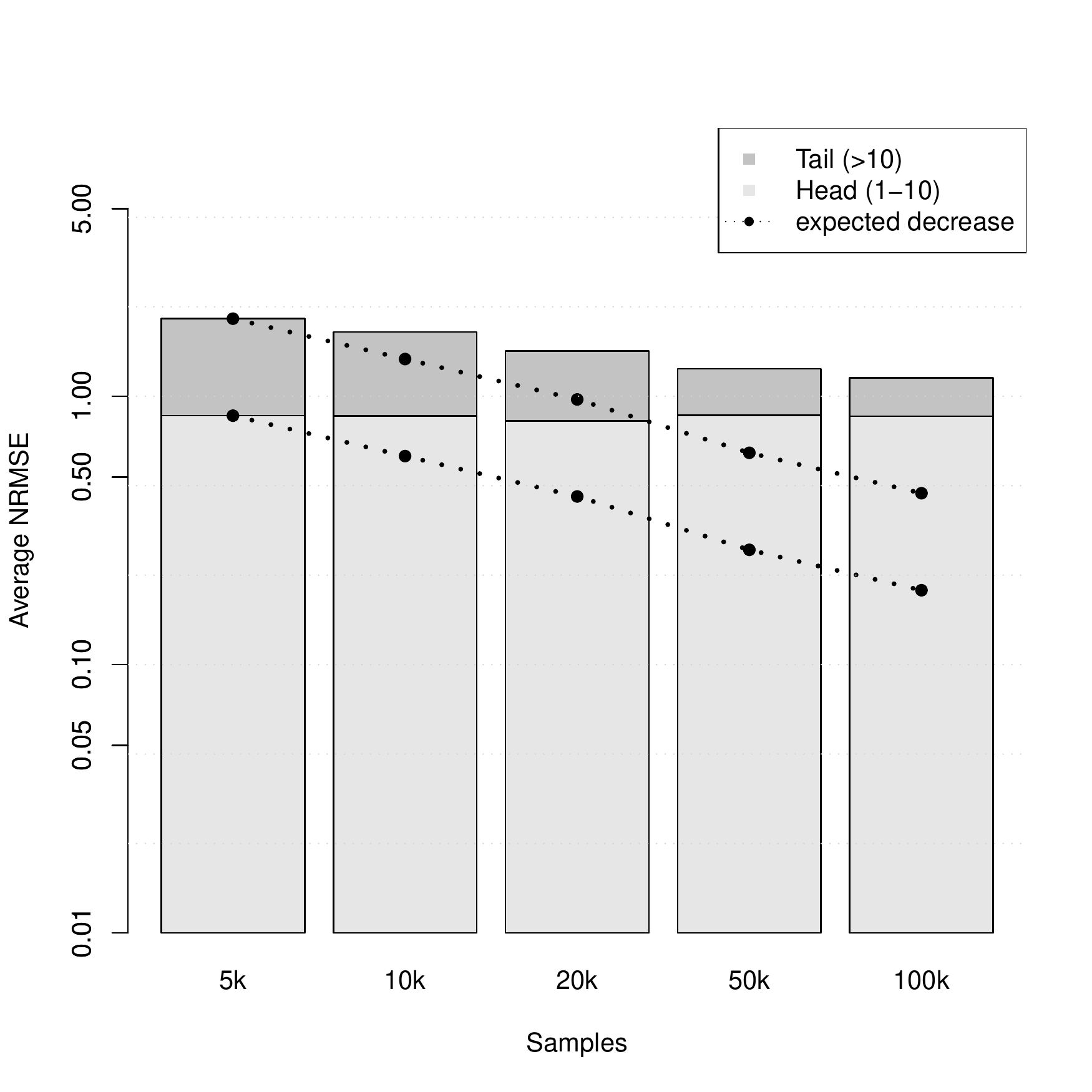}
\label{fig:nmse_smallp}}}

\centerline{
\subfloat[]{\includegraphics[width
=2.1in]{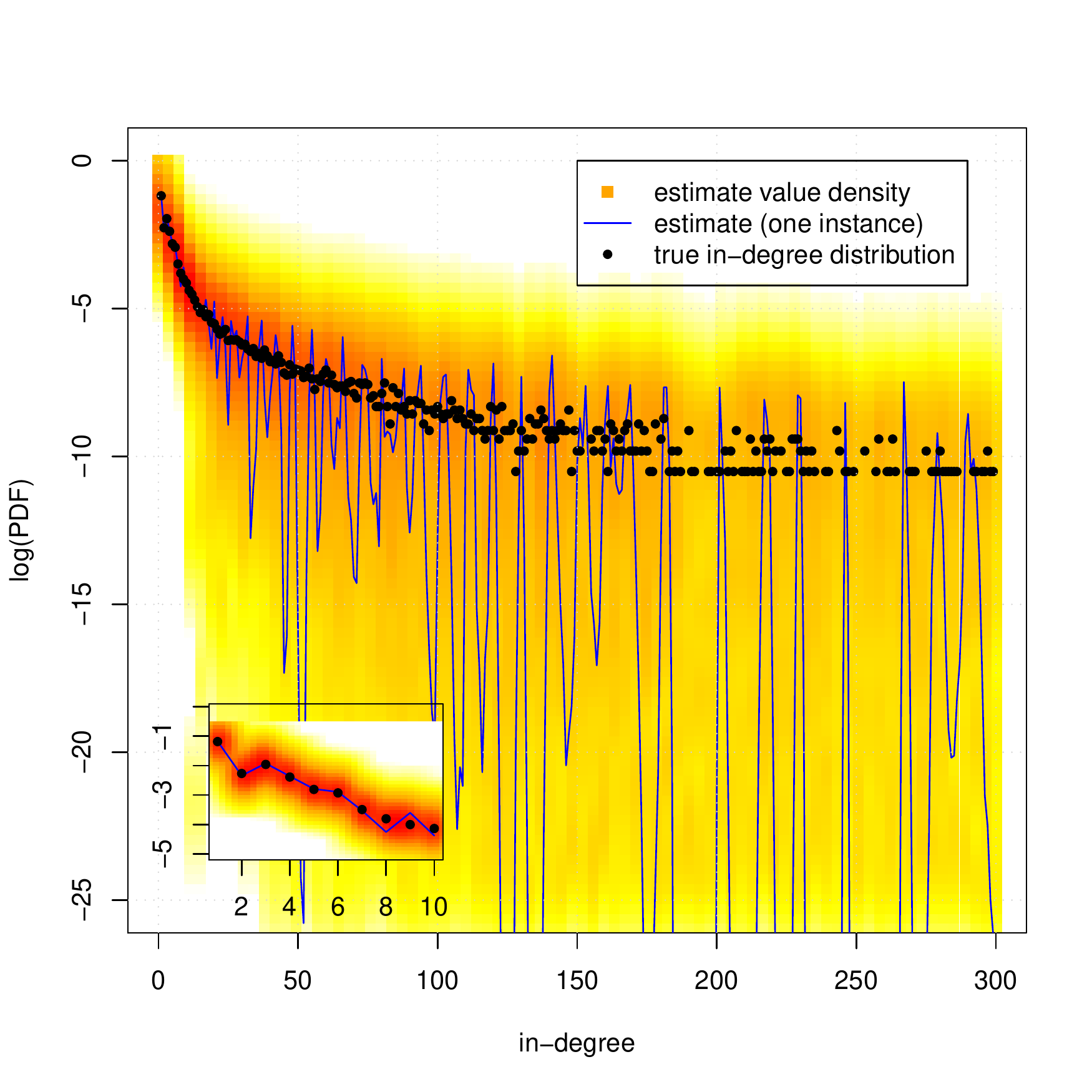}
\label{fig:heatmap_090_20k}}
\subfloat[]{\includegraphics[width
=2.1in]{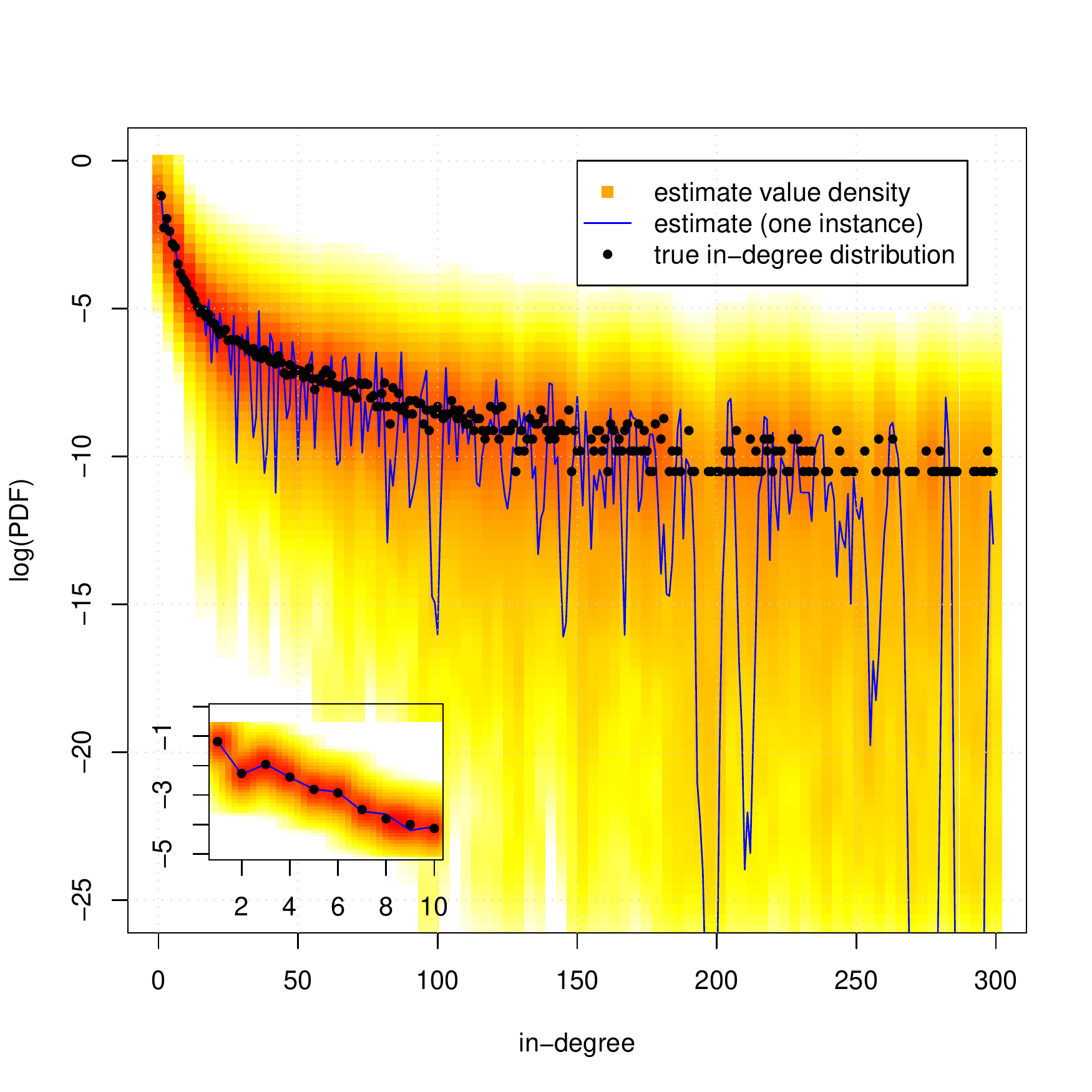}
\label{fig:heatmap_090_100k}}
\subfloat[]{\includegraphics[width
=2.1in]{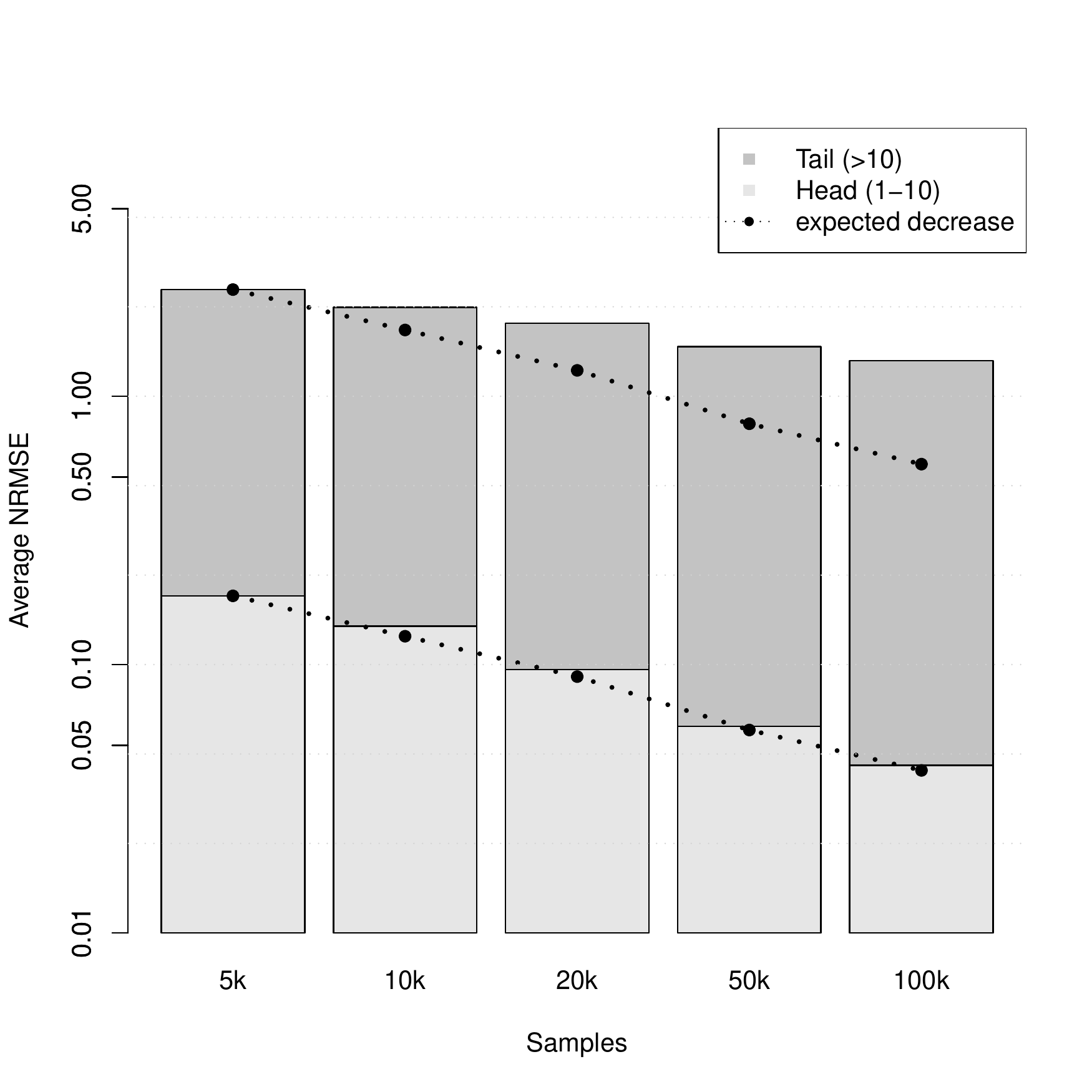}
\label{fig:nmse_largep}}
}
\caption{ The first row (a-c) shows the results for $p=0.25$, while the second row
    (d-e) shows the corresponding plots for $p=0.90$.
(a-b,d-e) True degree distribution, one example of estimate and
   heat map indicating the ocurrence rates of the estimate values for $N=10\times 10^3$
   samples
   (first column)
   and $N=50\times 10^3$ samples (second column), respectively. The red color in
   the heat map indicates high density of estimated values and
yellow (white) indicates low (no) density of estimated values. A subplot shows a zoom-in for the first
 degrees.
(c,f) Average NRMSE of the head and the
 tail of the distribution for $N \in \{1,5,10,20,100\}\times 10^3$.
 Dashed line shows how the error should vary with the number of samples.
 In (c) we have the {\bf typical behavior of wrong estimates}.
Increasing the number of samples does not improve the
quality of estimates. On the other hand (f) shows the {\bf typical behavior of correct estimates}.
Here increasing the number of samples yields lower estimation errors of the head.
}
\label{fig:experiments}
\end{figure}

Collecting a fraction of the exchanged messages means sampling network edges. 
Disregarding edge weights, assume the directed edges are independently sampled with probability $p$.
Henceforth, each person with more than one observed incoming email shall be called a sample.
Figure~\ref{fig:heatmap_025_10k} depicts the quality of the in-degree estimator in~\eqref{eq:mle} (see Section~\ref{sec:results} for the derivation) with $p=0.25$, leading to $N = 10^4$ sampled individuals. The black dots indicate the true in-degree distribution, the
blue curve shows a typical estimate, and the heat map indicates the
density of estimated values across 100 runs, where red indicates high density and
yellow (white) indicates low (no) density of estimated values. 
We observe from the blue curve that the estimated values can be orders of
magnitude away from the actual values and from the heat map we observe that the blue line is typical.

In what follows we illustrate the effects of varying the number of samples $N$
or changing the sample probability $p$ separately.
To vary $N$ while keeping $p$ fixed, we draw a node in-degree directly from the
in-degree distribution of this network and subsequently sample its edges. We
repeat this process until we obtain $N$ observed sets. This can be seen as sampling a
larger (smaller) network that has the same degree distribution.

We make two main observations:
\begin{enumerate}
\item {\bf Increasing the number of samples yields {\em no} reduction in estimation errors.}
This is an odd behavior.
We know from estimation theory that the error should decrease by $\sqrt{M}$ when the number
of samples is increased by a factor of $M$.
Figure~\ref{fig:heatmap_025_50k} shows the corresponding results for $N=50\times 10^3$.
We observe that the estimated fraction of nodes of each degree can still be very far from the actual values.

To make it clear that the accuracy gain from increasing the number of
samples is not in agreemeent with theory, we compute the estimate error obtained when we vary the
number of samples $N \in \{5,10,20,50,100\}\times 10^3$, for $p=0.25$.
The error is first measured in terms of the Normalized Root Mean Square Error
(NRMSE), which is defined as
$$
\textrm{NRMSE}(\hat{\theta_i})=\frac{\sqrt{E[(\hat{\theta_i}-\theta_i)^2]}}{\theta_i}.
$$
where $\hat{\theta}_i$ and $\theta_i$ are the estimated and true fraction of degree $i$ nodes,
respectively.
Then we take the average NRMSE from the head (degrees up to 10) and
the tail (degrees larger than 10) of the distribution separately.

Surprisingly, we observe in Figure~\ref{fig:nmse_smallp} that
there is almost no improvement in accuracy across different sample
sizes, even when we compare $5\times 10^3$ and $10^5$ samples. We also display in this
figure the expected reduction in the NRMSE for both head and tail by dashed lines.
It turns out that the error
does not decrease as we would expect. This raises the question of why, which we
address in Section~\ref{sec:results}.



\item {\bf For much larger values of $p$, the error starts to decrease with the number of samples.}
According to Theorem~\ref{thm:main} that we describe in Section~\ref{sec:results},
the difficulties experienced above arise due to the use of
small sampling probability ($p<0.5$) with heavy-tailed distributions,
ant not due to a lack of samples.
Hence we repeat the experiment using
$p=0.9$.
Figures~\ref{fig:heatmap_090_20k} and \ref{fig:heatmap_090_100k} show the
heat maps for $N=20\times 10^3$ and $N=10^5$. As opposed to what we previously saw,
increasing the number of samples makes the estimates closer to the true
in-degree distribution. The accuracy gain as a function of the
number of samples is shown in
Figure~\ref{fig:nmse_largep}. In fact, we observe that the NRMSE does decrease
as expected for the head of the distribution, but not for the tail.
Why are there two distinct behaviors, one for the head and one for the tail?
Why did it help to
increase the number of samples when estimating frequencies of small degrees
for $p=0.9$, as opposed to what we observed for $p=0.25$? Is it possible
to make the NRMSE of the tail to decrease as fast as the NRMSE of the head?
\end{enumerate}

In order to investigate the questions we pose here,
we study the Cram\'er-Rao Lower Bound (CRLB) of the set size
estimation problem. This give us a lower bound
on the estimation errors based on the amount of
information contained in the samples, measured
in terms of Fisher Information. Moreover, we
apply the CRLB to the estimation of the in-degree
distribution and average in-degree.

\section{Model} \label{sec:model}
Let $\cS_k$ be a nonempty set of elements, $k=1,\dots,m$, with $\cS_i \cap \cS_j = \emptyset$, $i,j =1,\ldots m$, $i\ne j$.
Let $S_k = \vert \cS_k \vert$ denote the size of the $k$-th set and assume set sizes are i.i.d.\ with distribution
$S_k \sim \btheta=(\theta_1,\dots,\theta_W)$, $W > 1$ $k \geq 1$. 
We assume $W$ finite ($W < \infty$).
The model breaks nodes (edges) into groups (sets) and our task in what follows is to characterize those groups from incomplete observation (sample) of these sets.
To illustrate the model, consider a directed graph; the set of incoming (outgoing) edges of a node $k$ is represented by $\cS_k$, 
$\btheta$ is the indegree (outdegree) distribution, and $W$ is the maximum indegree (outdegree).
Another straightforward example is representing IP traffic of a communications network, where $k$ is a TCP flow, $\cS_k$ is the set of TCP/IP packets that constitute flow $k$, and $W$ is the maximum observable flow size.

\subsection*{Sampling}
We observe (sample) elements of $\cS_k$, $k=1,\dots,m$, with probability $p$ -- a process also known as thinning.
Let $\alpha(\cS_k)$ be a random function that returns the number of observed elements of $\cS_k$ .
Elements are sampled independently (i.e., the sampling process is Bernoulli) and thus,
\[
  P[ \alpha(\cS_k) = j \vert \, S_k = i] = \begin{cases}
                                                       \binom{i}{j} p^j q^{i-j} \, , & j\geq 0, i > 1, i \geq j , \\
                                                       0, & \text{otherwise},
                                                       \end{cases}
\] 
where $q=1-p$.
We assume that when no elements of a set are observed, then the set as a whole
is not observed, i.e., $\cS_k$ is said to be observable if $\alpha(\cS_k) > 0$.
Thus, we denote $$\bbS = \{\alpha(\cS_k) : \alpha(\cS_k) > 0 \, , \, k=1,\ldots, m\}$$ the size of the observable set sizes.
Let $N = |\bbS|$ denote the number of observed sets.

\subsection*{Estimation}
We start by considering $p=1$, that is, all elements of all sets are observed.
The minimum variance estimator of $\theta_i$ is 
\[
 T_i^\prime(\cS_1,\ldots, \cS_m) = \sum_{k=1}^m \frac{{\bf 1}\{S_k = i\}}{N},  \]
where $N=m$. 
To measure the accuracy of the estimates we consider the mean squared error (MSE)  -- a.k.a.\ quadratic loss -- of the estimates
\[
  \MSE(T_i^\prime(\cS_1,\ldots,\cS_m)) = E[(T_i^\prime(\cS_1,\ldots,\cS_m) - \theta_i)^2] = \frac{\theta_i(1-\theta_i)}{ m } \leq \frac{1}{4 m} .\] 
Thus, for $p=1$ the estimation error decreases as $1/m$, recalling that $m$ is the number of sets.

Unfortunately, accurately estimating $\btheta$ when $p < 1$ is significantly more challenging.
Recall that a set $\cS_k$ is said to be observable if $\alpha(\cS_k) > 0$.
{\em We upfront assume that a unobservable sets cannot 
be used in the estimation process.}
This means that our estimator only has access to sets $\cS_k$ where $\alpha(\cS_k) > 0$.
Here we need another function $T_i$ that takes the observed set sizes $\bbS$ as inputs and outputs an {\bf unbiased} estimate $T_i(\bbS)$ of $\theta_i$, i.e., $E[T_i(\bbS)] = \theta_i$.
In what follows we focus on unbiased estimates; our discussion section (Section~\ref{sec:discussion}) extends our results to biased estimators.
The Mean Squared Error (MSE) of our estimator is
\[
  \MSE(T_i(\bbS)) =  E[(T_i(\bbS) - \theta_i)^2].
\]
The function $T_i$ that minimizes the MSE with respect to sets of size $i=1,\ldots,W$ is
\[
T^\star_i(\bbS) = \argmin_{T_i} \MSE(T_i(\bbS)) ,
\]
s.t.\ $E[T^\star_i(\bbS)] = \theta_i$.

\section{Results} \label{sec:results}
In this section we present and discuss our results.
\begin{theorem}\label{thm:main}
Let $\btheta=(\theta_1,\dots,\theta_W)$ be a distribution where $\exists i_0$
such that $\theta_i \leq 1/2$ for all $i > i_0$.
Recall that $N \leq m$ is the number of observed sets out of the total $m$ sets.
We show that, as $W \rightarrow \infty$, for $N$ sufficiently large any unbiased
estimator $T_i(\bbS)$, $i\geq 1$ is such that:
\begin{enumerate}
\item When $\theta_W$ decreases faster than exponentially in $W$, i.e.,
$-\log \theta_W = \omega (W)$, $\MSE(T_i(\bbS)) = O(1/N)$ for $0 < p <1$.
\item When $\theta_W$ decreases exponentially in $W$, i.e., $\log \theta_W = W\log a + o(W)$ as for some $0 < a < 1$,
\begin{enumerate}
\item $\log[\MSE(T_i(\bbS))] = \Omega (W/\log N)$, if $p < a/(a+1)$,
\item  $\MSE(T_i(\bbS)) = \Omega(W^{2 i+1}/N)$, if $p = a/(a+1)$, 
\item  $\MSE(T_i(\bbS)) = O(1/N)$, if $p > a/(a+1)$.
\end{enumerate}
\item When $\theta_W$ decreases more slowly than exponential, i.e., $-\log \theta_W = o(W)$,
\begin{enumerate}
\item $\log[\MSE(T_i(\bbS))] = \Omega (W/\log N)$, if $p < 1/2$,
\item $\MSE(T_i(\bbS)) = O(1/N)$, if $p \geq 1/2$; more precisely, 
\begin{enumerate}
\item $\MSE(T_i(\bbS)) = \omega (1/N)$, if $p = 1/2$ and $\sum_{j=1}^W j^{2i}
\theta_j = \omega (1)$,
\item $\MSE(T_i(\bbS)) = O(1/N)$, if either $p > 1/2$ or $p = 1/2$ and
$\sum_{j=1}^W j^{2i} \theta_j = O(1)$.
\end{enumerate}
\end{enumerate}
\end{enumerate}
\end{theorem}

\begin{theorem} \label{thm:average}
The bounds on the estimation error of the average set size are
analogous to the set size distribution bounds.
\end{theorem}

In what follows we explain how we sketch out the proof of Theorems~\ref{thm:main} and~\ref{thm:average} and describe their implications.

\subsection{Lower Bound on Estimation Errors}
 In this section we derive a lower bound on the Mean Squared Error (MSE) of
 $T_i(\bbS)$, $i=1,\dots,W$.  For this we use the Cram\'er-Rao (CR) lower
 bound of $T_i(\bbS)$, which gives the smallest MSE that any unbiased estimator
 $T_i$ can achieve.

Recall that a set is observable only if one or more of its elements are observable.
The probability that a (random) set $\cS$ is observed and has $j$ elements is defined as
\begin{equation}\label{eq:bji}
  b_{ji}(p) \equiv P\left[\alpha(\cS)=j \given \alpha(\cS) >0, |\cS|=i \right]= \frac{\binom{i}{j} p^j q^{i-j}}{1 - q^i} , \quad \text{if }0 < j \leq i \leq W ,
\end{equation}
and $b_{ji}(p) = 0$ otherwise, where $q= 1-p$.
Let $d_j(\btheta,p)$ denote the fraction of observed sets with exactly $j$ observed elements.
From~(\ref{eq:bji}) we have, $j=1,\ldots,W$,
\begin{align}
d_j(\btheta,p) &= P[\alpha(\cS)=j | \, |\cS|>0] \nonumber\\
& = \sum_{i=j}^{W}P[\alpha(\cS)=j| \alpha(\cS)>0, |\cS|=i]P[|\cS|=i |\,\alpha(\cS)>0]\nonumber\\
& = \sum_{i=j}^{W}b_{ji}(p)\phi_i(\btheta).
\label{eq:dj}
\end{align}
where 
\begin{equation} \label{eq:theta2phi}
    \phi_i(\btheta) = P[|\cS| = i \given \alpha(\cS) > 0] =
        \frac{\theta_{i}(1-q^{i})}{\sum_{k=1}^{W}\theta_{k}(1-q^{k})},
\end{equation}
is the distribution of the set sizes of the observed sets.
Or, in matrix notation,
$$
d(\btheta,p) = B(p) \bphi(\btheta),
$$
where $d(\btheta,p)=(d_1(\btheta,p),\dots,d_W(\btheta,p))^\textsf{T}$ and $B(p) = [b_{ji}(p)], j,i = 1,\dots,W$.
To illustrate the distribution $d(\btheta,p)$ in our model, note that for a random observed set $\cS$,
\[
  \alpha(\cS) \sim d(\btheta,p),
\]
%
with likelihood function
\begin{equation}
    f(j | \btheta) \equiv P[\alpha(\cS)=j \given \btheta] = (B(p) \bphi(\btheta))_j =d_j(\bphi(\btheta),p).
\label{eq:mle}
\end{equation}
In what follows for simplicity we denote $d_j(\btheta,p)$ as $d_j(\btheta)$, $j=1,\ldots, W$.

Recall that we are interested in functions $T_i(\bbS)$ that take as input the observed subset sizes $\bbS$ and outputs an unbiased estimate $T_i(\bbS)$ of $\theta_i$, $i=1,\ldots,W$.
Moreover, we want these estimates to be accurate, i.e., $\text{MSE}(T_i(\bbS))$ must be low in respect to $\theta_i$.
Otherwise, the estimate is of little use to the practitioner for set sizes of
interest, as illustrated in Figure~\ref{fig:experiments}.

Thus, it is important to find attainable lower bounds of $\text{MSE}(T_i(\bbS))$.
The Cram\'er-Rao Theorem states that the MSE of
any unbiased estimator $T$ is lower bounded by the inverse of the
Fisher information matrix divided by the number of independent samples $N$, provided some weak regularity conditions hold~\cite[Chapter 2]{vanTrees}, i.e.,
\begin{equation} \label{eq:CR}
\text{MSE}(T_i(\bbS)) \equiv E[(T_i(\bbS) - \theta_i)^2] \geq  \frac{\left((J^{(\btheta)}(p))^{-1} \right)_{ii}}{N}, 1\leq i \leq W.
\end{equation}
where $(J^{(\btheta)}(p))^{-1}$ is the inverse of the Fisher information matrix of a {\em single observed set} defined using the likelihood function~\eqref{eq:mle} as
\begin{equation}
(J^{(\btheta)}(p))_{i,k} \equiv \sum_{j=1}^W \frac{\partial \ln f(j \given \btheta)}{\partial \theta_i} \frac{\partial \ln f(j\given \btheta)}{\partial \theta_k} d_j(\bphi(\btheta)) = \sum_{j=1}^W \frac{\partial d_j(\bphi(\btheta))}{\partial \theta_i} \frac{\partial d_j(\bphi(\btheta))}{\partial \theta_k}  \frac{1}{d_j(\bphi(\btheta))},
\label{eq:Jtheta}
\end{equation}
given  $\sum_{i=1}^W \theta_i = 1$.

The lower bound in~(\ref{eq:CR}) is known in the literature as the Cram\'er-Rao lower bound or {\em CRLB} for short.
{Let $T_i^*(\bbS)$ be an unbiased estimator, $i=1,\ldots$.
We say $T_i^*(\bbS)$ is asymptotically efficient if
$\textrm{MSE}(T_i^*(\bbS))$ approaches the Cram\'er-Rao lower bound in~\eqref{eq:CR} as $N \to \infty$. We show in Appendix~\ref{sec:mleefficiency} that the Maximum Likelihood Estimator is asymptotically efficient on the set size estimation.
The implication of having an efficient estimator is that the lower bounds provided in this paper are tight for $N$ sufficiently large.
In what follows we represent $J^{(\btheta)}(p)$ as
$J^{(\btheta)}$ for simplicity}.

\subsection{Obtaining the CRLB} \label{sec:obtainingCRLB}
In what follows we derive closed-form lower bounds for the MSE of any unbiased estimator $T$, as a function of the original set size distribution $\btheta$, the sampling probability $p$, and the number of observed sets $N$, where we ignore the constraint $\sum_{i=1}^W \theta_i = 1$. 
Deriving a closed-form solution for the inverse of $J^{(\btheta)}$ is no easy task as matrix $J^{(\btheta)}$ is a function of $\partial f(j|\btheta)/\partial \theta_i$, $i=1,\ldots,W$, which makes $J^{(\btheta)}$ a non-linear function of $\btheta$.
However, observe that the likelihood function $f^\star(j | \bphi) \equiv P[\alpha(\cS) = j | \bphi ]$ (where $\cS$ is a random observed set) is linear with respect to $\bphi$
\begin{equation}
    f^\star(j | \bphi) \equiv (B \bphi)_j =d_j(\bphi).
\label{eq:like_phi}
\end{equation}
It is worth noting that $f^\star(j|\bphi(\btheta)) = f(j|\btheta)$.
The Fisher information matrix with respect to $\bphi$ is defined as  
$J^{(\bphi)}=[J_{i,k}^{(\bphi)}], i,k = 1,\dots,W$, where
\begin{equation}\label{eq:Jphi}
J_{i,k}^{(\bphi)} \equiv \sum_{j=1}^W \frac{\partial d_j(\bphi)}{\partial \phi_i} \frac{\partial d_j(\bphi)}{\partial \phi_k}  \frac{1}{d_j(\bphi)},
\end{equation}
given $\sum_{i=1}^W \phi_i = 1$; 
and because $d_j(\bphi)$ is linear in $\bphi$, combining~\eqref{eq:like_phi} and~\eqref{eq:Jphi} yields
\begin{equation} \label{eq:Jphi_matrix}
(J^{(\bphi)})^{-1} = B(p)^{-1}\textrm{diag}(B(p)\bphi)^{-1}(B(p)^{-1})\transpose  - \bphi \bphi\transpose.
\end{equation}
Here the term $\bphi \bphi\transpose$ corresponds to the accuracy gain obtained by considering the constraint $\sum_{i=1}^W \phi_i = 1$ (see Tune and Darryl~\cite{TuneIMC08} for more details and Gorman and Hero~\cite{GormanHeroToIT90} for the general formula on adding equality constraints to the CRLB).
Quantitatively we can safely ignore the constant term $\bphi \bphi\transpose$ as we are interested in the behavior of $(J^{(\bphi)})^{-1}$ as a function of $W$ and the elements of $\bphi \bphi\transpose$ are typically small.
All that is left to do is to find a relationship between $(J^{(\bphi)})^{-1}$ and $(J^{(\btheta)})^{-1}$.

We now obtain $(J^{(\btheta)})^{-1}$ from $(J^{(\bphi)})^{-1}$ through a multi-variate extension of the single variable chain rule.
As $f^\star(j | \bphi(\btheta)) = f(j | \btheta)$ the chain rule yields
\begin{equation*}
\frac{\partial f(j | \btheta)}{\partial \theta_i} = \frac{\partial f^\prime(\phi_j(\theta_i)) }{\partial \theta_i} =
\frac{\partial f^\prime(\phi_j)}{\partial \phi_j}\cdot\frac{\partial \phi_j(\btheta)}{\partial \theta_i}, \quad \forall i,j.
\end{equation*}
Using the Jacobian  $\nabla H = [h_{ik}]$, $h_{ik}={\partial \theta_k(\bphi)}/{\partial \phi_i}$ with $\theta_k(\bphi)$ as given in~\eqref{eq:theta2phi},  we arrive at the equivalent multivariate rule~\cite[pp.\ 83]{vanTrees} to express $(J^{(\btheta)})^{-1}$ as 
\begin{equation} \label{eq:inverseJtheta}
(J^{(\btheta)})^{-1} = \nabla H (J^{(\bphi)})^{-1}\nabla H\transpose.
\end{equation}
Using~\eqref{eq:Jphi_matrix} -- detailed derivation relegated to the Appendices -- we find:
\begin{equation} \label{eq:Jij}
[(J^{(\bphi)})^{-1}]_{ij} =\sum_{k=\max(i,j)}^{W} \left(\frac{q}{p} \right)^{2k} \binom{k}{j} \binom{k}{i} (-1)^{-i-j}(q^{-i}-1)(q^{-j}-1) d_k(\btheta) .
\end{equation}
Substituting~\eqref{eq:Jij} into~\eqref{eq:inverseJtheta} -- and through a
variety of algebraic manipulations detailed in the Appendices -- yields
\begin{align} \label{eq:Jtheta_body}
[(J^{(\btheta)})^{-1}]_{ii} & =  \frac{1}{\eta^{2}}\Bigg(\underbrace{\frac{1}{(1-q^{i})^{2}}[(J^{(\bphi)})^{-1}]_{ii}}_{A_1(i)}+\underbrace{\theta_{i}^{2}\sum_{j=1}^{W}\sum_{k=1}^{W}\frac{[(J^{(\bphi)})^{-1}]_{kj}}{(1-q^{k})(1-q^{j})}}_{A_2(i)} \nonumber\\
 &   - \underbrace{2\theta_{i}\sum_{j=1}^{W}\frac{[(J^{(\bphi)})^{-1}]_{ij}}{(1-q^{j})(1-q^{i})}}_{A_3(i)} \Bigg),
 \end{align}
where $\eta=\sum_{j=1}^{W}\phi_{j}(\btheta)/(1-q^{j})$. 
Note that term $A_1(i)$ of~\eqref{eq:Jtheta_body} is proportional to the CRLB of $\bphi$, $[(J^{(\bphi)})^{-1}]_{ii}$ but terms $A_2(i)$ and $A_3(i)$ are more involved.
Through a series of algebraic manipulations of terms $A_1$, $A_2$, and $A_3$,
        all detailed in the Appendices, we see that $(A_1(i) + A_2(i) - A_3(i))$ grows as a function of $(1-p)/p$ and $W$, yielding the relation
\begin{equation} \label{eq:TOmega}
\text{MSE}(T_i(\bbS)) = \Omega \left( \frac{\sum_{j=1}^W \left( \frac{1-p}{p} \right)^j \theta_j }{N}\right), \quad i=1,\ldots,W,
\end{equation}
where the number of observed sets $N$ is large but constant in respect to $W$.

The result in~\eqref{eq:TOmega} is very powerful as it gives simple  estimation
error lower bounds as a function of the sampling probability $p$ and the
original set size distribution $\btheta$. A close look at~\eqref{eq:TOmega}
reveals -- a detailed exposition is presented in the Appendices -- that when $((1-p)/p)^i \theta_i = \Omega(i^{-1})$ for all $i > i^\star$, $i^\star \ll W$, then the sum in~\eqref{eq:TOmega} grows at least as fast as the a harmonic series, which grows as $\log W$.
On the other hand, we see in the Appendices that when $((1-p)/p)^i \theta_i = O(i^{-\beta})$, $\beta > 1$, then the sum in~\eqref{eq:TOmega} converges to a constant, more precisely, it grows no faster than a Riemman zeta function with parameter $\beta$, $\zeta(\beta)$.

Thus, for a given $\btheta$ with  $W \gg 1$ the CRLB suffers from an interesting sharp threshold related to the sampling probability $p$.
If $p$ is below this threshold no estimator $T_i$ of $\theta_i$ ,$i=1,\ldots,W$, is able to achieve accurate estimates of $\theta_i$.
Below such $p$ threshold, and as long as the number of sampled sets, $N$, is large enough, there exists estimators $T_i(\bbS)$ ,$i=1,\ldots,W$, that can achieve accurate estimates.
To be more specific, we look at the threshold behavior of $p$ by breaking down $\btheta$ into three broad classes of distributions:
\begin{enumerate}
\item
If $\theta_W$ decreases faster than exponentially in $W$ there is no threshold
behavior of $p$. This is because if  $-\log \theta_W = \omega (W)$, then there
exists a constant $a < 1$ such that $((1-p)/p)^j \theta_j < a^j$, $j =
1,2,\ldots$. Hence, the sum in~\eqref{eq:TOmega} converges to a constant for any
$p > 0$, yielding $\text{MSE}(T_i(\bbS)) = \Omega(1/N)$, for $0 < p <1$.
Detailed arguments are presented in the Appendices.
\item
If $\log \theta_W = W\log a + o(W)$ then if $p \leq a/(a+1)$ yields $((1-p)/p)^j
\theta_j = a^{-j}\theta_j  = \Omega(1)$, $\forall j$. Hence, the sum
in~\eqref{eq:TOmega} diverges with $W$. On the other hand, if $p > a/(a+1)$ the
sum in~\eqref{eq:TOmega} converges to a constant. Detailed arguments are
presented in the Appendices.
\item
Finally, if $\theta_W$ decreases more slowly than exponential then if $p = 1/2 - \epsilon$, $\epsilon \geq 0$, yields $((1-p)/p)^j > (1+\epsilon/2)^j$, $\forall j$. Hence,  because $\theta_j$ decreases more slowly than an exponential, the sum in~\eqref{eq:TOmega} diverges with $W$. 
If $p \geq 1/2$ the lower bound in~\eqref{eq:TOmega} converges to a constant.
Detailed arguments are presented in the Appendices.

\end{enumerate}

To illustrate our results, we compute the MSE lower bounds in~\eqref{eq:Jtheta_body} 
where $\btheta$ is the Enron in-degree distribution truncated at different values of $W$. 
More precisely, we take the in-degree distribution of the Enron dataset (discussed in
Section~\ref{sec:motivation}) and truncate the maximum degree to $W$ by
accumulating in $W$ all the probability mass previously corresponding to
degrees greater than $W$.
The Enron in-degree distribution is a (truncated) heavier-than-exponential distribution.

\begin{figure}[!t]
\centerline{
\subfloat[$p=0.25$]{\includegraphics[width
=2.5in]{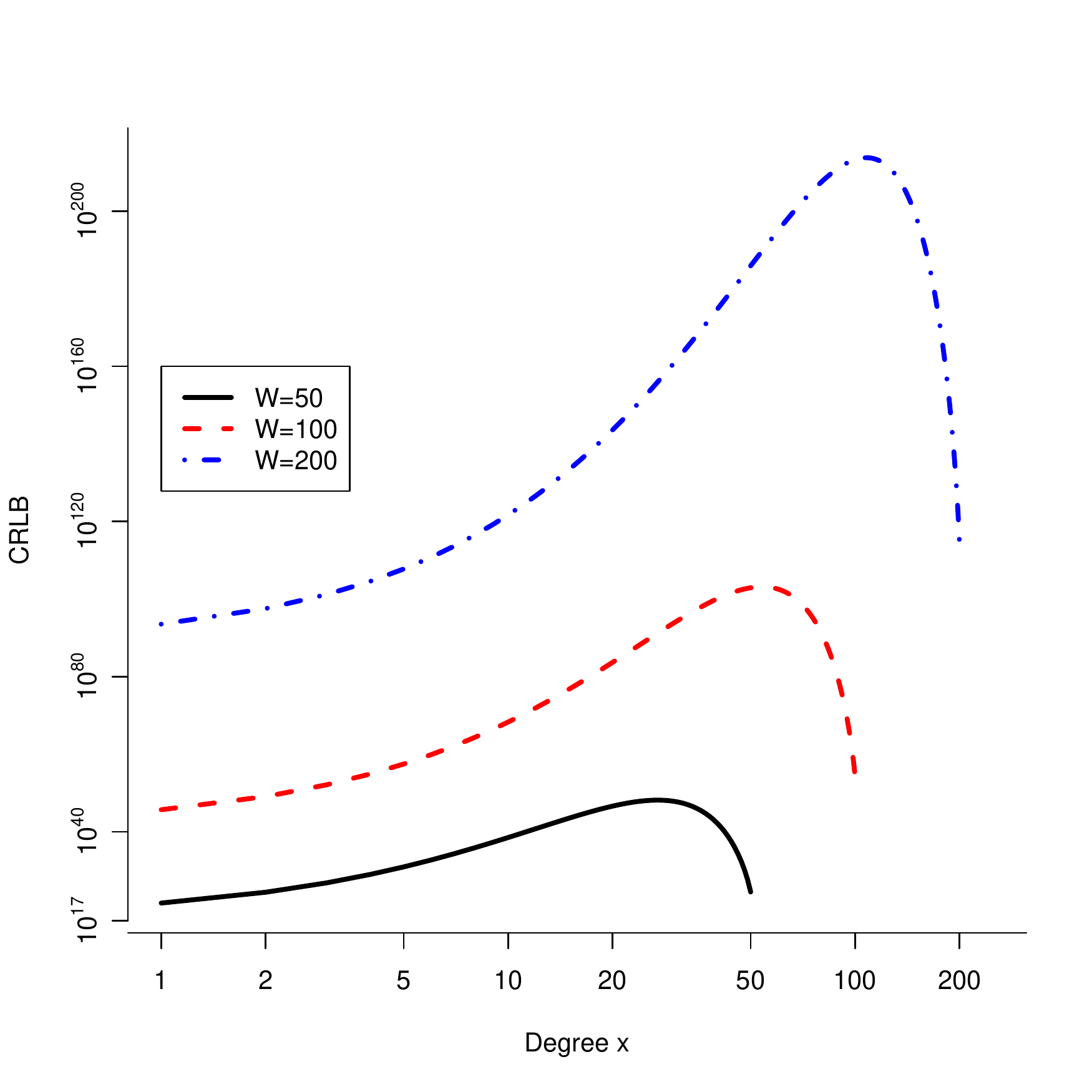}
\label{fig:crlb_p025}}
\subfloat[$p=0.90$]{\includegraphics[width
=2.5in]{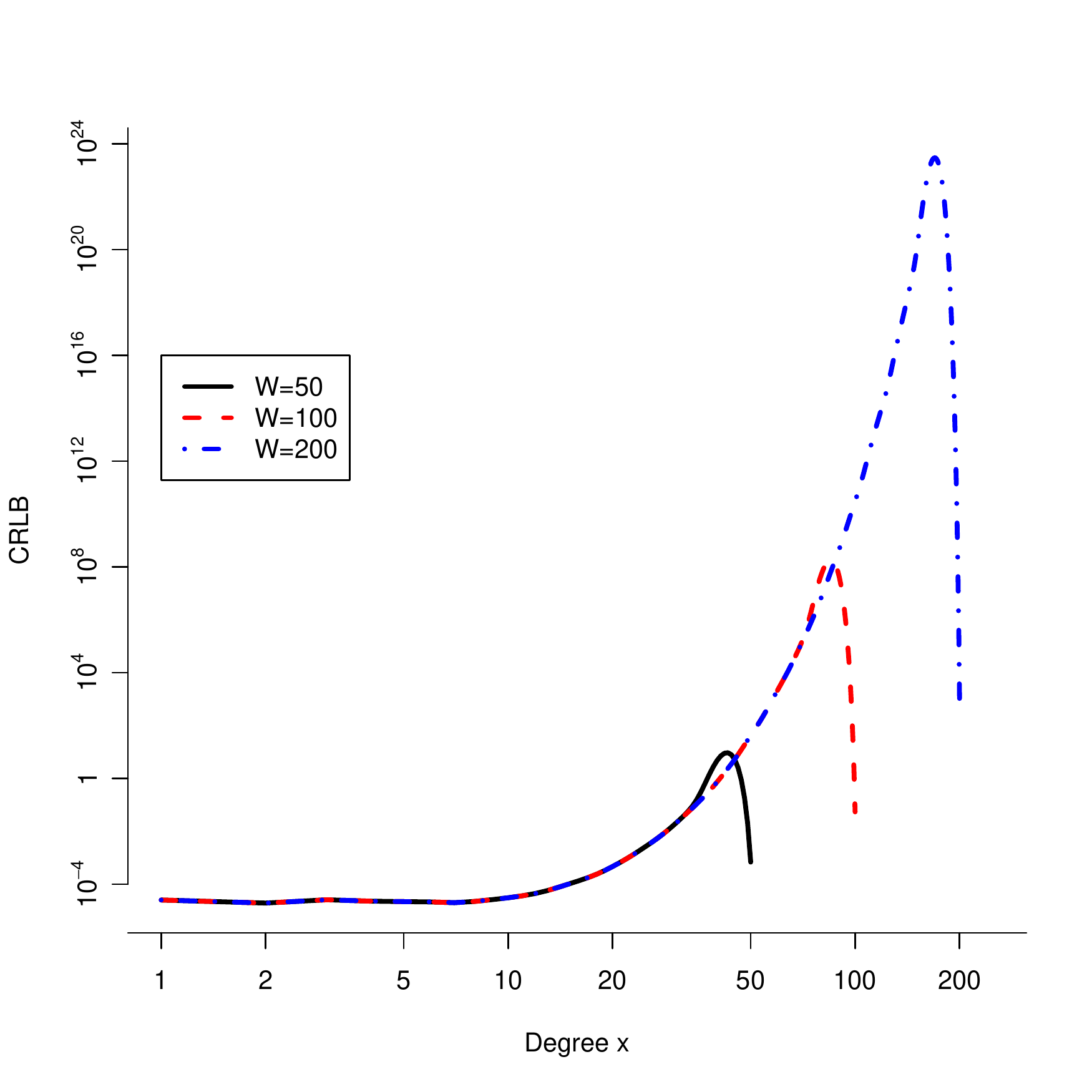}
\label{fig:crlb_p090}}
}
\caption{CRLB of the in-degree distribution of the Enron dataset for $N=10^4$
    samples.}
\label{fig:bounds}
\end{figure}

Figures~\ref{fig:bounds}a and~\ref{fig:bounds}b show the MSE lower bounds for
$p \in \{0.25,0.90\}$, respectively. We observe that for $p=0.25$ (Figure~\ref{fig:bounds}(a))
the MSE lower bound grows with $W$ even for small degrees, as predicted by Theorem~\ref{thm:main}.
While, for $p=0.9$ (Figure~\ref{fig:bounds}(b)) the MSE lower bound behaves
(mostly) independent of $W$, also as predicted by Theorem~\ref{thm:main}.These
results corroborate to explain the simulations results in
Section~\ref{sec:motivation}.

Other metrics besides the set size distribution are of interest.
In what follows we observe that, surprisingly, the accuracy of the average set size follows similar lower bounds of set size distribution estimators $T_i$, $i=1,\ldots,W$.
We then analyze the accuracy of entropy estimates.

\section{Accuracy of Estimated Averages}

In this section we focus on the accuracy of the average set size.
\subsection{Average set size}
The average set size is $m_\btheta = \sum_{j=1}^W j \theta_j$, or, alternatively, in matrix form
\[
   m_\btheta = [1,\ldots,W] \btheta\transpose .
\]
Let 
\[
   \frac{\nabla M}{\nabla\theta}  = \left[\frac{\partial m_\btheta}{\partial \theta_1}, \cdots , \frac{\partial m_\btheta}{\partial \theta_W} \right] = [1,\ldots,W].
\]
%
%
Let $m(\bbS)$ be an unbiased estimate of the average set size.
Using a similar argument used to obtain~\eqref{eq:inverseJtheta} (see Appendices) yields
\begin{align}
MSE(m(\bbS)) & \geq  \frac{\nabla M}{\nabla\theta} (J^{(\btheta)})^{-1} \frac{\nabla M}{\nabla\theta}\transpose\nonumber \\
   & = \frac{\nabla M}{\nabla\theta}\left(\frac{\nabla
        H}{\nabla\phi}(J^{(\bphi)})^{-1}\frac{\nabla H}{\nabla\phi}\transpose \right)\frac{\nabla M}{\nabla\theta}\transpose \nonumber \\
 & =  \left(\frac{\nabla M}{\nabla\theta}\frac{\nabla H}{\nabla\phi}\right)(J^{(\bphi)})^{-1}\left(\frac{\nabla M}{\nabla\theta}\frac{\nabla H}{\nabla\phi}\right)\transpose .\label{eq:MH}
\end{align}
Note that
\begin{align}
\left[\frac{\nabla M}{\nabla\theta}\frac{\nabla H}{\nabla\phi}\right]_{k} & = \sum_{i=1}^{W}ih_{ik}\nonumber \\
 & =  \sum_{{i=1\atop i\neq k}}^{W}i\left(-\frac{\theta_{i}}{\eta(1-q^{k})}\right)+k\left(\frac{1-\theta_{k}}{\eta(1-q^{k})}\right)\nonumber \\
 & =  \frac{1}{\eta(1-q^{k})}\left(k-\sum_{i=1}^{W}i\theta_{i}\right)\nonumber \\
 & =  \frac{k-m_{\theta}}{\eta(1-q^{k})},\label{eq:mh1}
 \end{align}
 where again $\eta=\sum_{j=1}^{W}\phi_{j}(\btheta)/(1-q^{j})$.
Substituting~\eqref{eq:mh1} into~\eqref{eq:MH} yields
\begin{align*}
\textrm{MSE}(m(\bbS)) & \geq 
\frac{1}{N}\sum_{i=1}^{W}\sum_{j=1}^{W}\left(\frac{j-m_{\theta}}{\eta(1-q^{j})}\right)[(J^{(\bphi)})^{-1}]_{ji}\left(\frac{i-m_{\theta}}{\eta(1-q^{i})}\right)\\
 & = \frac{1}{N} \frac{1}{\eta^{2}}\Bigg(\sum_{i=1}^{W}\sum_{j=1}^{W}\frac{ij[(J^{(\bphi)})^{-1}]_{ji}}{(1-q^{j})(1-q^{i})}+
 m_{\theta}^{2}\sum_{i=1}^{W}\sum_{j=1}^{W}\frac{[(J^{(\bphi)})^{-1}]_{ji}}{(1-q^{j})(1-q^{i})}-\\
 & 
 2m_{\theta}\sum_{i=1}^{W}\sum_{j=1}^{W}\frac{j[(J^{(\bphi)})^{-1}]_{ji}}{(1-q^{i})(1-q^{j})}\Bigg)  \\
 &  = \frac{1}{N} \frac{1}{\eta^{2}}\Bigg( \eta(\sum_{i=1}^{W}i^{2}\theta_{i}+\frac{q}{p}m_{\theta}) + \frac{m_{\theta}^{2}}{\theta_{j}^{2}}A_2(i) - 2m_{\theta}\eta(m_{\theta}+\frac{q}{p}\theta_{1})\Bigg),
 \end{align*}
 with $A_2(i)$ as given in~\eqref{eq:Jtheta_body}.  Detailed derivations are
 found in the Appendices.
 A closer look at $A_2(i)$ reveals
 \begin{equation} \label{eq:mOmega}
A_{2}(i)=\frac{1}{N}\theta_{i}^{2}\left(1+\eta\left(\sum_{j=1}^{W}q^{j}\theta_{j}+\sum_{j=1}^{W}\left(\frac{1-p}{p}\right)^{j}\theta_{j}\right)\right) = \Omega\left( \frac{\sum_{j=1}^{W}\left(\frac{1-p}{p}\right)^{j}\theta_{j}}{N} \right).
\end{equation}
Note that the lower bound of $m(\bbS)$ in~\eqref{eq:mOmega} is the same as the
lower bound of $T_i(\bbS)$, $i=1,\ldots,W$, in~\eqref{eq:TOmega}. Hence, a theorem in the lines of Theorem~\ref{thm:main} can be stated for $m(\bbS)$:
\begin{theorem}\label{thm:m}
Let $\btheta=(\theta_1,\dots,\theta_W)$ be a distribution where $\exists i_0$
such that $\theta_i \leq 1/2$ for all $i > i_0$.
Recall that $N \leq m$ is the number of observed sets out of the total $m$ sets.
We show that, as $W \rightarrow \infty$, for $N$ sufficiently large any unbiased estimator of the estimated mean of $\btheta$, $m(\bbS)$, must obey the following properties:
\begin{enumerate}
\item When $\theta_W$ decreases faster than exponentially in $W$, i.e., $-\log
\theta_W = \omega (W)$, $\MSE(m(\bbS)) = O(1/N)$ for $0 < p <1$.
\item When $\theta_W$ decreases exponentially in $W$, i.e., $\log \theta_W = W\log a + o(W)$ as for some $0 < a < 1$,
\begin{enumerate}
\item $\log[\MSE(m(\bbS))] = \Omega (W/\log N)$, if $p < a/(a+1)$,
\item  $\MSE(m(\bbS)) = \Omega(W/N)$, if $p = a/(a+1)$, 
\item  $\MSE(m(\bbS)) = O(1/N)$, if $p > a/(a+1)$.
\end{enumerate}
\item When $\theta_W$ decreases more slowly than exponential, i.e., $-\log \theta_W = o(W)$,
\begin{enumerate}
\item $\log[\MSE(m(\bbS))] = \Omega (W/\log N)$, if $p < 1/2$,
\item $\MSE(m(\bbS)) = O(1/N)$, if $p \geq 1/2$; more precisely, 
\begin{enumerate}
\item $\MSE(m(\bbS)) = \omega (1/N)$, if $p = 1/2$ and $\sum_{j=1}^W j^{2} \theta_j = \omega (1)$,
\item $\MSE(m(\bbS)) = O(1/N)$, if either $p > 1/2$ or $p = 1/2$ and $\sum_{j=1}^W j^{2} \theta_j = O(1)$.
\end{enumerate}
\end{enumerate}
\end{enumerate}
\end{theorem}

%
%
%
%
%
%

Theorem~\ref{thm:m} states that estimating the average set size is in the same order of hardness as estimating the entire set size distribution.

It is interesting, though, to verify if the same property holds in the case of
the average size of the observed sets, i.e., the
    average set size in respect to $\bphi$, $$m_\bphi =  \sum_{j=1}^W j \phi_j .$$
In what follows we show that the difficulty in estimating $m_\bphi$ is a
function of $W$ and is affected only by the first and second moments of $\bphi$, that is, as long as $m_\bphi$ and
\[
m^{(2)}_\bphi = \sum_{j=1}^W j^2 \phi_j
\]
are finite, $m_\bphi$ can be accurately estimated if enough samples, $N$, are collected.


Let $\hat m_\phi(\bbS)$ denote an unbiased estimate of $m_\phi$ and let
\[
\textrm{MSE}(\hat m_\phi(\bbS)) = E[(\hat m_\phi(\bbS) - m_\phi)^2]
\]
denote the MSE of $\hat m_\phi(\bbS)$.
After applying a variety of algebraic manipulations detailed in the Appendices we arrive at the following inequality
\begin{align*}\label{eq:msemphi}
\textrm{MSE}(\hat m_\phi) &\geq \frac{(1,\dots,W) (J^{(\phi)})^{-1}(1,\dots,W)^\texttt{T} - m_{\phi}^2}{N} \\
&=\sum_{k=1}^{W} \sum_{i=1}^k \sum_{j=1}^k  i j  \binom{k}{j} \binom{k}{i}
\frac{(-q)^{2k-i-j}}{p^{2k}}(1-q^i)(1-q^j) d_k(\bphi) \nonumber \\
&= \left( \sum_{i=1}^W \frac{i ( p i + q^{i+1} - 2q^i +q)\phi_i }{p(1-q^i)}  - m_\phi^2  \right)/N.
\end{align*}
More interestingly, we show that
\begin{equation}
\hat{m}^\star_\phi(\bbS) = \frac{\sum_{s \in \bbS} s}{N p}+\left(1-\frac{1}{p}\right) \frac{\sum_{s \in \bbS} {\bf 1}_{s=1}}{N},
\end{equation}
is an unbiased efficient (minimum variance) estimator of $m_\phi$, yielding
\[
\textrm{MSE}(\hat m^\star_\phi(\bbS)) =  \left( \sum_{i=1}^W \frac{i ( p i + q^{i+1} - 2q^i +q)\phi_i }{p(1-q^i)}  - m_\phi^2  \right)/N.
\]
Alternatively we can rewrite the above as
\[
\textrm{MSE}(\hat m_\phi) = O\left(\frac{m^{(2)}_\phi - m_\phi^2}{N}\right) .\]


Hence, $\MSE(\hat m_\phi)$ is lower bounded by the variance of the observed set sizes.
A simple explanation for this behavior is likely found in the inspection paradox.
Even if we know the sizes of the sampled sets, the mere fact that the set is sampled means that it probably has a higher than average size, as the probability that a set of size $i$ is sampled is $1-(1-p)^i$. 
Larger variance in the set sizes means larger biases towards sampling larger sets, which in turn makes it harder to unbias these samples.

\section{Discussion}\label{sec:discussion}

%
%
We divide this section in three parts.
Section~\ref{sec:estinit} considers the initialization of estimation procedures.
Section~\ref{sec:dpi} shows that no clever way to process the data $\bbS$ exists that would allow an estimator to violate the bounds provided in Section~\ref{sec:results}.
Finally, Section~\ref{sec:estbias} shows that our results can be extended to encompass biased and Bayesian estimators.

\subsection{Initialization of Estimation Procedures}\label{sec:estinit}
As previously stated, eq.~(\ref{eq:mle}) can be used to derive a maximum
likelihood estimator (MLE) for $\btheta$. From
the MLE one could either use a constrained non-linear optimization method to
maximize the likelihood function directly or use the Expectation-Maximization
(EM) algorithm to write an iterative estimation procedure. In the latter case,
the procedure consists of an initialization step followed by a loop of two
steps known as the E-step and M-step. We discuss two issues that arise when
EM is used to estimate the set size distribution.

In EM, the solution to which the algorithm converges to depends on the initial
guess. Therefore, in order to have an unbiased estimate, one must choose a point
uniformly at random from the space of possible values. 
Although it may seem reasonable to choose values for each $\theta_i$ uniformly in $[0,1]$ and
then normalize them, it turns out that this does not yield uniformly
distributed initial guesses. \fixme{One way to correctly generate the initial guess is to
draw from the Dirichlet
distribution with $W$ parameters $\mathbf{\alpha}=(1,\ldots, 1)$, since the Dirichlet PDF at
point $\btheta$ is proportional to $ \prod_{i=1}^W \theta_i^{\alpha_i-1}$.}

Nevertheless, such an initialization combined with the other two steps of EM will
give us estimates $\hat{\theta}_i \in [0,1]$ hence producing biased estimates.  Therefore, it is possible that EM achieves
an MSE not in agreement with the CRLB we derived previously. This is
the case when the number of samples $N$ is small and, consequently, the
diagonal of $G$ has relatively large values (possibly greater than 1). On
the other hand, for large $N$, the number of observed sets with size $i$
will converge to a Normal distribution with mean $\theta_i$ and small
variance. For small enough variance, restricting $\theta_i$ to be between
0 and 1 does not affect the final estimate significantly and thus the CRLB
accurately bounds the MSE.

\subsection{An Application of the Data Processing Inequality}\label{sec:dpi}
The data processing inequality~\cite{Zamir98} states that no function of the data may increase the amount of Fisher information already contained in the data.
Thus, the bounds in Theorems~\ref{thm:main} and~\ref{thm:m} remain unchanged regardless of how the data is pre-processed, no matter how clever the pre-processing approach is. This, of course, encompasses any type of noise filters or machine learning methods.

\subsection{Impact on Different Types of Estimators: Bayesian, Frequentist, Biased and Unbiased} \label{sec:estbias}
To extend our results beyond unbiased estimators we explain the connection between Fisher information, the Cram\'er-Rao bound and biased estimators.
We also extend our results to Bayesian estimators (including maximum a posteriori estimators).
\subsubsection{Extension to Biased Estimators}
Let $h(\theta_i) = E[T_i(\bbS)] - \theta_i$ be the estimator bias. Then (see for instance Ben-Haim and Eldar~\cite{ConstrainedBiasFI09})
\[
\MSE(T_i(\bbS)) \geq \left(1+\frac{\partial b(\theta_i)}{\partial \theta_i} \right)^2 [(J^{(\btheta)})^{-1}]_{ii},
\]
assuming $\partial b(\theta_i)/\partial \theta_i$ exists.
Note if the bias derivative satisfies $-2 < \partial b(\theta_i)/\partial
\theta_i < 0$, then the biased estimator has lower MSE than any unbiased
estimator. However, we believe it is unlikely that a large value of
$[(J^{(\btheta)})^{-1}]_{ii}$ (as large as $10^{160}$ as seen in
        Section~\ref{sec:obtainingCRLB} for the Enron e-mail network) can be compensated by a biased estimator.
\subsubsection{Extension to Bayesian Estimators}
Let $\btheta$ now be a random variable with prior distribution $\pi_\btheta$.
A Bayesian estimator adds $\pi_\btheta$ as extra information to the estimation problem.
The Fisher information of the prior is 
\[
 J^{(p)}_{ij} = E\left[ \frac{\partial \ln \pi_\btheta}{\partial \theta_i} \frac{\partial \ln \pi_\btheta}{\partial \theta_j} \right].
\]
The Fisher information obtained exclusively by the data is $J^{(\btheta)}$ presented in~\eqref{eq:Jtheta}.
And the total Fisher information {\em prior + data} is~\cite[pp.\ 84]{vanTrees}
\[
J^{(t)} = J^{(p)} + J^{(\btheta)}.
\]
The Cram\'er-Rao bound of a Bayesian estimator $W_i(\bbS)$ of $\theta_i$ with prior $\pi_\btheta$ yields~\cite[pp.\ 85]{vanTrees}
\[
  \MSE(W_i(\bbS))  \geq (J^{(t)})^{-1} = (J^{(p)} + J^{(\btheta)})^{-1},
\]
and thus, if the data contains \fixme{little Fisher information} then a decrease
in the MSE is due to
the information contained in the prior $\pi_\btheta$.

\section{Conclusions \& Related Work} \label{sec:conclusion}
In this paper we give explicit expressions of MSE lower bounds of unbiased estimators of the distribution of set sizes $\theta$ and the average set size $m_\btheta$ with sampling probability $p$.
We show that the estimation error of $\btheta$ grows at least exponentially in $W$, when $\log \theta_W = W\log a + o(W)$ as $W\rightarrow \infty$ for some $0 < a < 1$, and $p <a/(a+1)$, or when $\log \theta_W = o(W)$ as $W\rightarrow \infty$ and $p <1/2$, which indicates that there unbiased estimators of some distributions $\btheta$ are too inaccurate to be useful for practitioners. Moreover we show that unbiased estimates of $m_\theta$ suffer from similar problems. 

Not much prior work exists in the literature.  Hohn and Veitch \cite{Horn03} first
observed that using a sampling probability of $p < 1/2$ poses problems in the
context of two specific estimators for the flow size distribution when the
distribution obeys a power law.  In particular, they showed that their
estimators are asymptotically unbiased with decreasing error as the number of
flow samples increases when $p \geq 1/2$ but not when $p < 1/2$.  Our work shows that
this is a fundamental result of set size distribution estimation and not
specific to any one or two estimators. Ribeiro et al. \cite{RTTB06} was the first to
introduce the use of Fisher information as a design tool for flow size
estimation.  Experiments reported in that paper suggested that there is little
information when p is small and showed how this information can be significantly
increased with the addition of other data taken from packet headers. Last, Tune
and Veitch \cite{TuneIMC08} applied Fisher information to compare packet sampling with flow
sampling. In the process of doing so, they obtained a variety of useful Fisher
information inverse identities, which we rely on in this work.

\section{Acknowledgments}
This research was sponsored by the NSF under CNS-1065133, ARO under MURI
W911NF-08-1-0233, and the U.S. Army Research Laboratory under Cooperative Agreement
W911NF-09-2-0053. The views and conclusions contained in this document are
those of the authors and should not be interpreted as representing the official policies,
either expressed or implied of the NSF, ARO, ARL , or the U.S. Government. The U.S.
Government is authorized to reproduce and distribute reprints for Government purposes
notwithstanding any copyright notation hereon.

\bibliographystyle{plain}

\techreport{
\appendices

\section{Set size distribution proofs}

Let $B(p)=[b_{ji}(p)], j,i=1,\dots,W$ be a matrix whose elements are given by
\begin{equation}
  b_{ji}(p) \equiv P\left[\alpha(\cS)=j\given \alpha(\cS) >0, |\cS|=i \right]= \frac{\binom{i}{j} p^j q^{i-j}}{1 - q^i} , \quad \text{if }0 < j \leq i ,
\end{equation}
and $b_{ij}(p)=0$ otherwise, where $q=1-p$.

Lemma~\ref{lemma:Bstar} shows a closed formula for the inverse of $B(p)$.
\begin{lemma}\label{lemma:Bstar}
$B(p)^{-1}=[b_{ji}^\star(p)]$ ($i,j=1,\ldots, W$), where
\[
b_{ji}^\star(p) =
\begin{cases}
           \binom{i}{j} p^{-i}(-q)^{i-j}(1-q^j)& i \geq j \\
           0 & i<j.
\end{cases}
\]
\end{lemma}
\begin{pf}
Let $B(p)^{-1}=[b_{ji}^\star(p)]$ with $b_{ji}^\star(p)$ defined above.  We first
show that $Y=B(p)B(p)^{-1}$ is an identity matrix.
Consider element $(j,i)$ of $Y$:
\begin{equation}\label{eq:y}
y_{ji}=\sum_{l=1}^{W} b_{jl}(p) b_{li}^\star(p) \, .
\end{equation}
We have three cases: $j > i$, $j = i$, and $j < i$.\\
{\em Case 1, $j>i$}: eq.~(\ref{eq:y}) yields $y_{ji} = 0$ since $b_{jl}(p)=0$,
   $\forall l \leq i$ and $b_{li}^\star(p)=0$, $\forall l > i$.

\noindent{\em Case 2, $j=i$}:  Here $b_{jl}(p) b_{lj}^\star(p)=0$, $\forall l \neq j$ and (\ref{eq:y}) yields
\[
y_{jj}= \frac{p^j}{1-q^j}\cdot p^{-j} (1-q^j)=1 \, .
\]
\noindent {\em Case 3, $j<i$}:
 eq. (\ref{eq:y}) yields
\begin{equation*}
\begin{split}
y_{ji}&=\sum_{l=j}^{i} (-1)^{i-l} p^{j-i} q^{i-j} \binom{l}{j} \binom{i}{l}\\
&= p^{j-i}q^{i-j} \sum_{l=j}^{i} (-1)^{i-l} \binom{i}{j} \binom{i-j}{l-j}\\
&= p^{j-i}q^{i-j} \binom{i}{j} \sum_{l=j}^{i} (-1)^{i-l} \binom{i-j}{l-j}\\
&= p^{j-i}q^{i-j} \binom{i}{j} (1-1)^{i-j}\\
&=0\\
\end{split}
\end{equation*}
Thus, $y_{jj} = 1$, $\forall j$ and $y_{ji} = 0$, $\forall j \neq i$, which concludes our proof.
\end{pf}

Lemma~\ref{lemma:Bstar} directly yields the inverse of the Fisher information
matrix $J^{(\phi)}$ of a single observed set, as seen in the following lemma.

\begin{lemma}\label{lemma:inversematrix}
$(J^{(\phi)})^{-1}=[[(J^{(\phi)})^{-1}]_{ij}]$ $(i,j=1,2,\ldots, W)$, where
\begin{equation}
[(J^{(\phi)})^{-1}]_{ij}=\sum_{k=\max(i,j)}^{W} \left(\frac{q}{p}\right)^{2k}
\binom{k}{j} \binom{k}{i} (-1)^{-i-j}
    (q^{-i}-1)(q^{-j}-1) d_k(\btheta)
\end{equation}
\end{lemma}
\begin{pf}
Denote $R^{(\bphi)}(p)=[R^{(\bphi)}_{ji}(p)]=B^{-1}(p)\textrm{diag}(B(p)\bphi)^{-1}$, where
$R^{(\bphi)}_{ji}(p)=b_{ji}^\star(p) d_i(\bphi)$.
Based on Lemma~\ref{lemma:Bstar} and eq.~\eqref{eq:dj}, we have
\baq  \label{eq:R}
  R^{(\bphi)}_{ji}(p) = \left\{ \barr{ll}
  \binom{i}{j} p^{-i}(-q)^{i-j}(1-q^j) d_i(\bphi), & i \geq j, \\
  0, & i < j.
  \earr \right.
\eaq

Since $J^{(\phi)}=R^{(\bphi)}(p)(B(p)^{-1})\transpose$, $[(J^{(\phi)})^{-1}]_{ji}$ is computed as the following equation based on Lemma~\ref{lemma:Bstar} and eq. \eqref{eq:R}
\baqm
[(J^{(\phi)})^{-1}]_{ji} &=& \sum_{k=1}^{W} R^{(\bphi)}_{jk}(p) b_{ik}^\star(p) \\
&=&\sum_{k=\max(i,j)}^{W} \frac{\binom{k}{j} \binom{k}{i}
    (-q)^{2k-i-j}(1-q^i)(1-q^j) d_k(\bphi)}{p^{2k}}\\
&=&\sum_{k=\max(i,j)}^{W} \left(\frac{q}{p}\right)^{2k} \binom{k}{j} \binom{k}{i}
    (-1)^{-i-j}(q^{-i}-1)(q^{-j}-1) d_k(\bphi)
\eaqm
\end{pf}


\begin{lemma}\label{lemma:Gjj}
$(J^{(\theta)})^{-1}=[[(J^{(\theta)})^{-1}]_{ij}]$ $(i,j=1,2,\ldots, W)$, where
\beq
[(J^{(\theta)})^{-1}]_{ii} =
\frac{1}{\eta^{2}}\Bigg(\frac{[(J^{(\phi)})^{-1}]_{ii}}{(1-q^{i})^{2}}+\theta_{i}^{2}\sum_{j=1}^{W}\sum_{k=1}^{W}\frac{[(J^{(\phi)})^{-1}]_{kj}}{(1-q^{k})(1-q^{j})}
        -2\theta_{i}\sum_{j=1}^{W}\frac{[(J^{(\phi)})^{-1}]_{ij}}{(1-q^{i})(1-q^{j})}\Bigg)\label{eq:g_jj}
\eeq
where $\eta=\sum_{i=1}^{W}\phi_{i}/(1-q^{i})$. 
\end{lemma}

\begin{pf}
The relationship between $(J^{(\theta)})^{-1}$ and $(J^{(\bphi)})^{-1}$ is given
by
\begin{equation}
(J^{(\theta)})^{-1} = \nabla H (J^{(\bphi)})^{-1} \nabla H\transpose, 
\end{equation}
where $\nabla H = [h_{ik}]$ with $h_{ik}=\partial \theta_k(\bphi)/\partial \phi_i$.
Hence
\[
h_{ik}=\begin{cases}
-\frac{\phi_{i}/(\eta(1-q^{i}))}{\eta(1-q^{k})} & i\neq k\\
\frac{1-\phi_{i}/(\eta(1-q^{i}))}{\eta(1-q^{i})} & i=k\end{cases}\]
where $\eta=\sum_{k=1}^{W}\phi_{k}/(1-q^{k})$ is a constant. Note
that from eq.~(\ref{eq:theta2phi}) we have $\theta_i = \phi_{i}/(\eta(1-q^{i}))$.
Therefore the diagonal elements of $(J^{(\theta)})^{-1}$ can be written as
\baq
[(J^{(\theta)})^{-1}]_{ii} & = &
\sum_{j=1}^{W}\sum_{k=1}^{W}h_{ik}[(J^{(\bphi)})^{-1}]_{kj}h_{ij}\transpose\nonumber\\
 & = & \sum_{{j=1\atop j\neq i}}^{W}\sum_{{k=1\atop k\neq
     i}}^{W}\left(-\frac{\theta_{i}}{\eta(1-q^{k})}\right)[(J^{(\bphi)})^{-1}]_{kj}\left(-\frac{\theta_{i}}{\eta(1-q^{j})}\right)+\nonumber\\
& & \sum_{{j=1\atop j\neq
    i}}^{W}\left(\frac{1-\theta_{i}}{\eta(1-q^{i})}\right)[(J^{(\bphi)})^{-1}]_{ij}\left(-\frac{\theta_{i}}{\eta(1-q^{j})}\right)+\nonumber\\
 &  & \sum_{{k=1\atop k\neq
     i}}^{W}\left(-\frac{\theta_{i}}{\eta(1-q^{k})}\right)[(J^{(\bphi)})^{-1}]_{ki}\left(\frac{1-\theta_{i}}{\eta(1-q^{i})}\right)+
\left(\frac{1-\theta_{i}}{\eta(1-q^{i})}\right)^{2}[(J^{(\phi)})^{-1}]_{ii}\nonumber\\
&=&\frac{1}{\eta^{2}}\Bigg(\frac{[(J^{(\phi)})^{-1}]_{ii}}{(1-q^{i})^{2}}+\theta_{i}^{2}\sum_{j=1}^{W}\sum_{k=1}^{W}\frac{[(J^{(\phi)})^{-1}]_{kj}}{(1-q^{k})(1-q^{j})}
        -2\theta_{i}\sum_{j=1}^{W}\frac{[(J^{(\phi)})^{-1}]_{ij}}{(1-q^{i})(1-q^{j})}\Bigg).
\eaq

\end{pf}

We split eq.~(\ref{eq:g_jj}) in three parts to carry out its analysis:
\beq
 [(J^{(\theta)})^{-1}]_{ii} =
 \frac{1}{\eta^{2}}\Bigg(\underbrace{\frac{[(J^{(\theta)})^{-1}]_{ii}}{(1-q^{i})^{2}}}_{A_{1}(i)}+
         \underbrace{\theta_{i}^{2}\sum_{j=1}^{W}\sum_{k=1}^{W}\frac{[(J^{(\theta)})^{-1}]_{kj}}{(1-q^{k})(1-q^{j})}}_{A_{2}(j)}-
         \underbrace{2\theta_{i}\sum_{j=1}^{W}\frac{[(J^{(\theta)})^{-1}]_{ij}}{(1-q^{i})(1-q^{j})}}_{A_{3}(i)}\Bigg).
 \eeq

\subsection{Analysis of $A_{1}(i)$}

Based on Lemma~\ref{lemma:inversematrix} and eq.~\eqref{eq:dj}, we have
\begin{lemma}\label{lemma:Jjj}
\begin{equation}
A_{1}(i)=\eta q^{-2i}\sum_{j=0}^{W-i}\binom{i+j}{i}q^{j+i}\theta_{j+i}g_{ij}.\label{eq:Jrelatedwithgji}
\end{equation}
where $\eta=\sum_{k=1}^W \phi_k/(1-q^k)$ and
$g_{ij}=\sum_{k=0}^j \bn{i+k}{i}\bn{j}{k}(q/p)^{k+i}$.
\end{lemma}
\begin{pf}
\baq
[(J^{(\phi)})^{-1}]_{ii} & = & \sum_{k=i}^W\left(\frac{q}{p}\right)^{2k}
\bn{k}{i}^2 (-1)^{-2i} (q^{-i}-1)^2 d_k(\bphi)  \nonumber \\
  & = & \sum_{k=i}^W \sum_{j=k}^W \left(\frac{q}{p}\right)^{2k}
\bn{k}{i}^2 (-1)^{-2i} (q^{-i}-1)^2 \frac{\bn{j}{k} p^k q^{j-k}\phi_j}{1 - q^j} \nonumber \\
  & = & (q^{-i}-1)^2 \sum_{j=i}^W \bn{j}{i}\frac{q^j\phi_j}{1-q^j}\sum_{k=i}^j \bn{k}{i}\bn{j-i}{k-i}(q/p)^k  \nonumber \\
  & = & (q^{-i}-1)^2 \sum_{j=0}^{W-i}
  \bn{i+j}{i}\frac{q^{i+j}\phi_{i+j} g_{ij}} {1-q^{i+j}}
\eaq
where $g_{ij}=\sum_{k=0}^j \bn{i+k}{i}\bn{j}{k}(q/p)^{i+k}$.

Since $\phi_{i}/(1-q^{i})=\theta_{i}\cdot \eta$, we can eq.~\eqref{eq:Jrelatedwithgji} as a function of $\theta$:

\[
[(J^{(\bphi)})^{-1}]_{ii}=\eta\left(q^{-i}-1\right)^{2}\sum_{j=0}^{W-i}\binom{i+j}{i}q^{i+j}\theta_{i+j}g_{ij}.\]
Therefore

\begin{equation}
A_{1}(i)=\eta q^{-2i}\sum_{j=0}^{W-i}\binom{i+j}{i}q^{i+j}\theta_{i+j}g_{ij}.\label{eq:t1}\end{equation}
\end{pf}

\begin{lemma} \label{lemma:general_bounds}
We have the following bounds for $A_1(i)$:
\beq
A_1(i)  <  C_i\sum_{k=0}^i c_{ik} \sum_{j=0}^{\infty} \bone{k\le
    j}(i+j)^{2i}
\bigl(\frac{q}{p}\bigr)^{i+j} \theta_{i+j} \label{eq:Jub}
\eeq
and
\beq \label{eq:Jlb}
A_1(i) > C_ic_{ii} \sum_{j=i(i-1)}^{W-i} j^{2i}
\bigl(\frac{q}{p}\bigr)^{i+j}\theta_{i+j}
\eeq
where
\[
C_i= \frac{\eta q^{-i}}{(i!)^2}
\]
and
\[ c_{ik} = \bn{i}{k}q^k \prod_{l=0}^{i-k-1} (i-l), \quad k=0,\ldots
,i;i=1,\ldots  W. \]
\end{lemma}
\begin{pf}
Since the $i$-th derivative of $(q/p)^{i+k}$ with respect to $q/p$, is
\[
\frac{\text{d}^i (q/p)^{i+k}}{\text{d}(q/p)^i}=\prod_{l=1}^i (k+l) (q/p)^k,
\]
we have the following equations for $g_{ij}$
\baqm
\lefteqn{ g_{ij} = \frac{1}{i!} \bigl( \frac{q}{p}\bigr)^i \sum_{k=0}^j
    \prod_{l=1}^i (k+l) \bn{j}{k}(q/p)^k }\\
 & = & \frac{1}{i!} \bigl( \frac{q}{p}\bigr)^i \sum_{k=0}^j \bn{j}{k}
 \frac{\text{d}^i (q/p)^{i+k}}{\text{d}(q/p)^i} \\
 & = & \frac{1}{i!} \bigl( \frac{q}{p}\bigr)^i \frac{\text{d}^i
     \Bigl(\sum_{k=0}^j \bn{j}{k} (q/p)^{i+k}\Bigr)}{\text{d}(q/p)^i} \\
 & = & \frac{1}{i!} \bigl( \frac{q}{p}\bigr)^i \frac{\text{d}^i \Bigl((q/p)^i (1
             + q/p)^j \Bigr)}{\text{d}(q/p)^i}.
\eaqm
Using a general form of the product rule~\cite[pp.\ 318]{Olver98} yields
\begin{equation}\label{eq:gji}
g_{ij} = \frac{1}{i!} \bigl( \frac{q}{p}\bigr)^i \sum_{k=0}^{\min\{i,j\}}
\bn{i}{k} \bigl( \frac{1}{p}\bigr)^{j-k} \prod_{l=0}^{k-1} (j-l)
    \bigl(\frac{q}{p}\bigr)^{k} \prod_{l=0}^{i-k-1} (i-l),
\end{equation}
where to simplify the expression we define $\prod_{l=0}^{-1} \cdots = 1$.

Substituting~\eqref{eq:gji} back into~\eqref{eq:t1}, we obtain the following
expression for $A_1(i)$
\begin{equation}\label{eq:JstarCloseForm}
A_1(i) = C_i \sum_{k=0}^i c_{ik} \sum_{j=0}^{W-i} \bone{k\le j}
\prod_{l=1}^{i}(j+l)\prod_{l=0}^{k-1}(j-l) (q/p)^{i+j}\theta_{i+j}
\end{equation}
where
\[
C_i= \frac{\eta q^{-i}}{(i!)^2}
\]
and
\[ c_{ik} = \bn{i}{k}q^k \prod_{l=0}^{i-k-1} (i-l), \quad k=0,\ldots
,i;i=1,\ldots , W. \]

We have the following upper bounds for $A_1(i)$,
\baq
A_1(i) & < &  C_i\sum_{k=0}^i c_{ik} \sum_{j=0}^{W-i} \bone{k\le j} (i+j)^{2i}
\bigl(\frac{q}{p}\bigr)^{i+j}\theta_{i+j}  \\
& < & C_i\sum_{k=0}^i c_{ik} \sum_{j=0}^{\infty} \bone{k\le j}(i+j)^{2i}
\bigl(\frac{q}{p}\bigr)^{i+j} \theta_{i+j}.
\eaq
A lower bound is obtained by noting that
\baqm
\prod_{l=1}^{i}(j+l)\prod_{l=0}^{k-1}(j-l) & > &  j^{i-k}
\prod_{l=1}^{k}(j+l)\prod_{l=1}^{k}(j-l+1)  \\
 & = &  j^{i-k} \prod_{l=1}^{k}(j^2 +j+l-l^2).
\eaqm
The latter is  greater than or equal to $j^{2i}$ whenever $j > i(i-1)$ yielding
\beq
A_1(i) > C_ic_{ii} \sum_{j=i(i-1)}^{W-i} j^{2i}
\bigl(\frac{q}{p}\bigr)^{i+j}\theta_{i+j}.
\eeq
\end{pf}


\subsection{Analysis of $A_{2}(i)$}

\begin{eqnarray}
\sum_{i=1}^{W}\sum_{j=1}^{W}\frac{[(J^{(\theta)})^{-1}]_{ij}}{(1-q^{i})(1-q^{j})}
 & = &
 \sum_{i=1}^{W}\sum_{j=1}^{W}\sum_{k=1}^{W}\frac{\binom{k}{j}\binom{k}{i}\left(\frac{q}{p}\right)^{2k}(-1)^{-j-i}(q^{-j}-1)(q^{-i}-1)d_{k}(\bphi)}{(1-q^{j})(1-q^{i})}\nonumber \\
 & = & \sum_{k=1}^{W}\left(\frac{q}{p}\right)^{2k}d_{k}(\bphi)\sum_{i=1}^{k}\sum_{j=1}^{k}\binom{k}{j}\binom{k}{i}(-q)^{-j-i}\nonumber \\
 & = & \sum_{k=1}^{W}\left(\frac{q}{p}\right)^{2k}d_{k}(\bphi)\left(\sum_{i=1}^{k}\binom{k}{i}(-q)^{-i}\right)^{2}\nonumber \\
 & = & \sum_{k=1}^{W}\left(\frac{q}{p}\right)^{2k}d_{k}(\bphi)\left(\left(-\frac{q}{p}\right)^{-k}-1\right)^{2}\qquad\textrm{using (\ref{eq:id1})}\nonumber \\
 & = & \sum_{k=1}^{W}d_{k}(\bphi)-2\sum_{k=1}^{W}\left(-\frac{q}{p}\right)^{k}d_{k}(\bphi)+\sum_{k=1}^{W}\left(\frac{q}{p}\right)^{2k}d_{k}(\bphi)\nonumber \\
 & = & 1-2\sum_{k=1}^{W}\left(-\frac{q}{p}\right)^{k}d_{k}(\bphi)+\sum_{k=1}^{W}\left(\frac{q}{p}\right)^{2k}d_{k}(\bphi).\label{eq:t2}\end{eqnarray}
First, note that

\begin{eqnarray}
\sum_{k=1}^{W}\left(-\frac{q}{p}\right)^{k}d_{k}(\bphi) & = &
\sum_{k=1}^{W}\left(-\frac{q}{p}\right)^{k}\sum_{j=1}^{W}\binom{j}{k}p^{k}q^{j-k}\theta_{j}\eta\nonumber \\
 & = & \eta\sum_{j=1}^{W}q^{j}\theta_{j}\sum_{k=1}^{j}\binom{j}{k}(-1)^{k}\nonumber \\
 & = & -\eta\sum_{j=1}^{W}q^{j}\theta_{j}.\qquad\textrm{using (\ref{eq:id3})}\label{eq:aux1}\end{eqnarray}
Also,

\begin{eqnarray}
\sum_{k=1}^{W}\left(\frac{q}{p}\right)^{2k}d_{k}(\bphi) & = &
\sum_{k=1}^{W}\left(\frac{q}{p}\right)^{2k}\sum_{j=1}^{W}\binom{j}{k}p^{k}q^{j-k}\theta_{j}\eta\nonumber \\
 & = & \eta\sum_{j=1}^{W}q^{j}\theta_{j}\sum_{k=1}^{j}\binom{j}{k}\left(\frac{q}{p}\right)^{k}\nonumber \\
 & = & \eta\sum_{j=1}^{W}q^{j}\theta_{j}\left(\left(\frac{1}{p}\right)^{j}-1\right)\qquad\textrm{using (\ref{eq:id2})}\nonumber \\
 & = &
 \eta\left(\sum_{j=1}^{W}\left(\frac{q}{p}\right)^{j}\theta_{j}-\sum_{j=1}^{W}q^{j}\theta_{j}\right).\label{eq:aux2}\end{eqnarray}
Replacing eqs. (\ref{eq:aux1}) and (\ref{eq:aux2}) into (\ref{eq:t2})
yields

\begin{eqnarray}
\sum_{i=1}^{W}\sum_{j=1}^{W}\frac{[(J^{(\theta)})^{-1}]_{ij}}{(1-q^{i})(1-q^{j})}
 & = &
 1+\eta\left(2\sum_{j=1}^{W}q^{j}\theta_{j}+\sum_{j=1}^{W}\left(\frac{q}{p}\right)^{j}\theta_{j}-\sum_{j=1}^{W}q^{j}\theta_{j}\right)\\
 & = &
 1+\eta\left(\sum_{j=1}^{W}q^{j}\theta_{j}+\sum_{j=1}^{W}\left(\frac{q}{p}\right)^{j}\theta_{j}\right).\label{eq:t2-final}\end{eqnarray}
Therefore,

\beq
A_{2}(i)=\theta_{i}^{2}\left(1+\eta\left(\sum_{j=1}^{W}q^{j}\theta_{j}+\sum_{j=1}^{W}\left(\frac{q}{p}\right)^{j}\theta_{j}\right)\right).
\label{eq:t2_exp}
\eeq

Note that $A_{2}(i)$ is positive and may diverge or not depending on the
summation $\sum_{j=1}^W \left(\frac{q}{p}\right)^j \theta_j$.

\subsection{Analysis of $A_{3}(i)$}

Note that

\begin{eqnarray}
\sum_{k=1}^{W}\binom{k}{i}\left(-\frac{q}{p}\right)^{k}d_{k}(\bphi)
 & = &
 \sum_{k=i}^{W}\binom{k}{i}\left(-\frac{q}{p}\right)^{k}\sum_{j=1}^{W}\binom{j}{k}p^{k}q^{j-k}\theta_{j}\eta\nonumber \\
 & = &
 \eta\sum_{k=i}^{W}(-1)^{k}\sum_{j=1}^{W}\binom{j}{i}\binom{j-i}{k-i}q^{j}\theta_{j}\nonumber \\
 & = &
 \eta\sum_{j=i}^{W}\binom{j}{i}q^{j}\theta_{j}\sum_{k=i}^{j}\binom{j-i}{k-i}(-1)^{k}\nonumber \\
 & = &
 (-1)^{i}\eta\sum_{j=i}^{W}\binom{j}{i}q^{j}\theta_{j}\sum_{k=0}^{j-i}\binom{j-i}{k}(-1)^{k}\nonumber \\
 & = & (-q)^{i}\eta\theta_{i}.\qquad\textrm{using (\ref{eq:id4})}\label{eq:aux3}\end{eqnarray}

We also have

\begin{eqnarray}
\sum_{k=1}^{W}\binom{k}{i}\left(\frac{q}{p}\right)^{2k}d_{k}(\bphi)
 & = &
 \sum_{k=1}^{W}\binom{k}{i}\left(\frac{q}{p}\right)^{2k}\sum_{j=1}^{W}\binom{j}{k}p^{k}q^{j-k}\theta_{j}\eta\nonumber \\
 & = &
 \eta\sum_{k=1}^{W}\left(\frac{q}{p}\right)^{k}\sum_{j=1}^{W}\binom{j}{i}\binom{j-i}{k-i}q^{j}\theta_{j}\nonumber \\
 & = &
 \eta\sum_{j=i}^{W}\binom{j}{i}q^{j}\theta_{j}\sum_{k=i}^{j}\binom{j-i}{k-i}\left(\frac{q}{p}\right)^{k}.\label{eq:aux4}\end{eqnarray}

From eq. (\ref{eq:aux3}) and (\ref{eq:aux4}), we have

\beq
\sum_{j=1}^{W}\frac{[(J^{(\theta)})^{-1}]_{ij}}{(1-q^{j})(1-q^{i})}=
\eta\theta_{i}-(-q)^{-i}\eta\sum_{j=i}^{W}\binom{j}{i}q^{j}\theta_{j}\sum_{k=i}^{j}\binom{j-i}{k-i}\left(\frac{q}{p}\right)^{k}\label{eq:aux5}\eeq

and hence,
\begin{equation} A_3(i) =
\underbrace{2\eta\theta_{i}^{2}}_{A_{3,1}(i)}-\underbrace{2\theta_{i}(-q)^{-i}\eta\sum_{j=0}^{W-i}
    \binom{i+j}{i}q^{i+j}\theta_{i+j}\sum_{k=0}^{j}\binom{j}{k}\left(\frac{q}{p}\right)^{k+i}}_{A_{3,2}(i)}.\label{eq:t3-final}\end{equation}

Since $A_{3,1}(i)$ is always finite, we only need to compare the magnitude
of $A_{1}(i)$ and $A_{3,2}(i)$. Since
$\sum_{k=0}^{j}\binom{j}{k}\left(\frac{q}{p}\right)^{k+i}<g_{ij}$,
we can bound $|A_{3,2}(i)|$ by

\[
|A_{3,2}(i)|\leq2\theta_{i}q{}^{-i}\eta\sum_{j=0}^{W-i}\binom{i+j}{i}q^{i+j}\theta_{i+j}g_{ij}.\]

Therefore

\[
A_{1}(i)-|A_{3,2}(i)|\geq(q^{-2i}-2\theta_{i}q{}^{-i})\eta\sum_{j=0}^{W-i}\binom{i+j}{i}q^{i+j}\theta_{i+j}g_{ij}.\]

The RHS of the previous inequation is positive when

\begin{eqnarray*}
q^{-2i} & \geq & 2\theta_{i}q{}^{-i}\\
\theta_{i} & \leq & \frac{1}{2q^{i}} < \frac{1}{2}.\end{eqnarray*}

{
Recall that we assumed that $\exists i_0$ such that $\theta_i \leq 1/2$
for all $i > i_0$. Thus by examining only $A_{1}(i)$ and $A_{2}(i)$ we can
determine whether $[(J^{(\theta)})^{-1}]_{ii}$ diverges or not for $i > i_0$.
}

\section{Proof of Theorem \ref{thm:main}.}

The lower bound of $\textrm{MSE}(T_i(\bbS))$, given by
$[(J^{(\theta)})^{-1}]_{ii}$, is described for each of the three possible cases
in Theorem~\ref{thm:main}. The corresponding proofs
are shown in what follows.

\textbf{ 1) When $\theta_{W}$ decreases faster than exponentially in
$W$.}

\begin{pf} 
Suppose that $\theta_{W}$ decreases faster than exponentially in
$W$. More precisely, assume that $-\log \theta_W = \omega (W)$.
It follows that  $\log (\theta_W / \theta_{W+1})
    \rightarrow \infty$ as $W\rightarrow \infty$.  Hence, for any $\epsilon >
    0$, there exists a $W_0(\epsilon )$ such that $\log (\theta_W/ \theta_{W+1}) >
    1/ \epsilon$ for $W>W_0(\epsilon )$. This implies $\theta_{W+1}/ \theta_W <
    e^{-1/\epsilon}$  for $W>W_0(\epsilon )$.  Given $p>0$, we can choose
    $\epsilon$ such that $q e^{-1/\epsilon} /p < 1$.   We now apply the ratio
    test for convergence of an infinite sum to each of the $i+1$ sums in the
    upper bound for $A_1(i)$ given by~\eqref{eq:Jub}.
\[
\frac{(W+i+1)^{2i} (q/p)^{W+i+1}\theta_{W+i+1}}{(W+i)^{2i} (q/p)^{W+i}\theta_{W+i}} <
\frac{(W+i+1)^{2i}}{(W+i)^{2i}}\frac{qe^{-1/\epsilon}}{p}
\]
for $W > W_0(\epsilon ) - i$ and the latter expression becomes less than one as
$W\rightarrow \infty$. Hence $A_1(i)= O(1)$ for $0<p<1$.

A similar argument can be used to show that $A_{2}(i) = O(1)$.
Hence, $[(J^{(\theta)})^{-1}]_{ii}=O(1)$ for $0<p<1$.
\end{pf}

\textbf{2) When $\theta_{W}$ decreases exponentially in $W$.}

\begin{pf}
Suppose that $\theta_{W}$ decreases exponentially in $W$. More precisely,
let $\log \theta_{W}=W \log a +o(W)$ for $0<a<1$.
Recall that $A_2(i)$ is positive. Therefore,
the logarithm of $[(J^{(\theta)})^{-1}]_{ii}$ in~\eqref{eq:g_jj} can be lower bounded as follows,
\beq
\log [(J^{(\theta)})^{-1}]_{ii} \ge \log A_1(i).
\eeq
In addition, the logarithm of $A_1(i)$ in~\eqref{eq:Jrelatedwithgji} can be
bounded by
\baqm  \log A_1(i) & \ge & W \log (q/p) + \log \theta_W + o(W) \\
 & = & W \log (qa/p)  + o(W)
\eaqm
where the latter equality follows from the hypothesis.  Now, if $qa/p >
1$, then $\log A_1(i) = \Omega (W)$, which implies
$\log [(J^{(\theta)})^{-1}]_{ii} = \Omega(W)$.  Note that $qa/p > 1$ iff $p<a/(a+1)$.

When $p=a/(a+1)$, then $qa/p = 1$.  Hence the lower bound of $A_1(i)$ given by \eqref{eq:Jlb}
is $\Omega (W^{2i+1})$. Hence, $[(J^{(\theta)})^{-1}]_{ii} =
\Omega(W^{2i+1})$.


Similarly to the proof for the case where $\theta_W$ decreases faster than
exponentially in $W$, we can use the ratio test for convergence of an infinite
sum to show that for $qa/p < 1$, $A_1(i) = O(1)$. Hence, it follows
that $[(J^{(\theta)})^{-1}]_{ii} = O(1)$ for $p>a/(a+1)$.

%
%
\end{pf}

\textbf{3) When $\theta_{W}$ decreases slower than exponentially in
$W$.}

\begin{pf}
Suppose that $\theta_{W}$ decreases slower than exponentially in $W$.
More precisely assume that
$-\log \theta_W = o(W)$.
The logarithm of $A_1(i)$ can be lower bounded as follows,
\baqm \log A_1(i) & \ge & W \log (q/p) + \log \theta_W + o(W) \\
 & = & W \log (q/p)  + o(W)
\eaqm
The latter equality follows from the hypothesis.  Now, if $q/p > 1$ (i.e., $p <
1/2$), then $\log
A_1(i) \ge \Omega (W)$, which implies
$\log  [(J^{(\theta)})^{-1}]_{ii} = \Omega(W)$.

When $p \ge 1/2$, it follows that $A_2(i) = O(1)$.
In particular if $p=1/2$ and $\sum_{j=1}^W j^{2i}\theta_j = \omega(1)$, we
can see from eq.~\eqref{eq:Jlb} that $A_1(i) = \omega(1)$ and in turn,
    $[(J^{(\phi)})^{-1}]_{ii} =\omega(1)$.

Note that for $p=1/2$ each of the $i+1$ sums in the upper bound for
$A_1(i)$ given by~\eqref{eq:Jub} is bounded by the ${2i}$-th moment of the
set size distribution. Hence, if $\sum_{j=1}^W j^{2i}\theta_j = O(1)$,
then $[(J^{(\theta)})^{-1}]_{ii} =O(1)$.

Finally, when $p > 1/2$, an argument similar to that used in the case where
$\theta_W$ decreases faster than exponentially yields $[(J^{(\theta)})^{-1}]_{ii} =O(1)$.
\end{pf}

\section{Simplified bounds}

It is worth noting that $A_2(i)$ gives us a lower bound on
$[(J^{(\theta)})^{-1}]_{ii}$, as $A_1(i)-A_3(i)>0$. Furthermore,
the convergence of $A_2(i)$ is given by the convergence of the sum $\sum_{j=1}^W
(q/p)^j \theta_j$. Therefore, we can write
\beq
[(J^{(\theta)})^{-1}]_{ii} = \Omega\left(\sum_{j=1}^W
\left(\frac{1-p}{p}\right)^j \theta_j\right).
\eeq
From that, we derive the following results.

\textbf{ 1) When $\theta_{W}$ decreases faster than exponentially in
$W$.}

By definition, for any $\epsilon>0$, there exists a $W_0(\epsilon)$ such that
$\log(\theta_W/\theta_{W+1}) > 1/\epsilon$. Given $p > 0$, we can choose
$\epsilon$ such that $qe^{-1/\epsilon}/p < 1$.  The ratio test for convergence
of an infinite sum reads
\beq
\frac{(q/p)^{j+1}\theta_{j+1}}{(q/p)^j\theta_j} < \frac{qe^{-1/\epsilon}}{p}
\eeq
Let $a=qe^{-1/\epsilon}/p$. Hence, there exists a $j^*$ such that
for all $j>j^\star$, $((1-p)/p)^j \theta_j <
a^j$, $j=1,2,\dots$. Therefore, the sum converges to a constant for any $0 < p < 1$,
    yielding $[(J^{(\theta)})^{-1}]_{ii} = O(1)$.

\textbf{2) When $\theta_{W}$ decreases exponentially in $W$.}

By definition, there exists $0 < a <1$ such that $\log \theta_W = W\log a +o(W)$.
When $p \leq a/(a+1)$ it follows that $((1-p)/p)^j\theta_j \ge a^{-j}\theta_j =
\Omega(1)$. Therefore, $[(J^{(\theta)})^{-1}]_{ii}= O(W)$. A tighter bound
can be obtained by taking into account $A_1(i)$, yielding $\log
[(J^{(\theta)})^{-1}]_{ii}= O(W)$ for $p < a/(a+1)$ and $
[(J^{(\theta)})^{-1}]_{ii}= O(W^{2i+1})$ for $p = a/(a+1)$. On the other
hand, for $p > a/(a+1)$, we have $((1-p)/p)^j\theta_j < a^{j}\theta_j =
O(1)$. Hence, $[(J^{(\theta)})^{-1}] = O(1)$.

\textbf{3) When $\theta_{W}$ decreases slower than exponentially in
$W$.}

When $p<1/2$, it follows that $(1-p)/p = a > 1$. In this case,
there exists a $j^\star$ such that for all $j>j^\star$,
$((1-p)/p)^j\theta_j = a^{j}\theta_j = \Omega(1)$. Hence,
$[(J^{(\theta)})^{-1}]_{ii} = O(W)$ for $p < 1/2$. Conversely,
when $p > 1/2$, $(1-p)/p = a < 1$. Hence, there exists a $j^\star$ such that for all $j>j^\star$,
$((1-p)/p)^j\theta_j = a^{j}\theta_j = O(1)$. Thus,
    $[(J^{(\theta)})^{-1}]_{ii}
= O(1)$ for $p > 1/2$. At last, for $p=1/2$, the summation is exactly 1, which
also implies $[(J^{(\theta)})^{-1}]_{ii} = O(1)$. In the latter case (i.e.,
$p=1/2$),
a tigher bound is obtained by taking $A_1(i)$ into account, which yields
$[(J^{(\theta)})^{-1}]_{ii} = \omega(1)$ if $\sum{j=1}^W j^{2i} \theta_j =
\omega(1)$ and 
$[(J^{(\theta)})^{-1}]_{ii} = O(1)$ if $\sum{j=1}^W j^{2i} \theta_j = O(1)$.

\section{Asymptotic Efficiency and Asymptotic Normality of the MLE $T_i^*(\bbS)$}
\label{sec:mleefficiency}

In this section we show that there exists a Maximum Likelihood Estimator (MLE)
$T_i^{(\phi)}(\bbS)$ of
$\phi_i$ that is asymptotic efficient (i.e., $\textrm{MSE}(T_i^*(\bbS)) =
[(J^{(\phi)})^{-1}]_{ii}$) and asymptotic normal. Since the Delta Method is an
exact approximation for the Normal distribution, it follows that there exists a
MLE $T_i^*(\bbS)$ of $\theta_i$  that is asymptotic efficient, which can be
obtained by applying the Delta Method to $T_i^{(\phi)}(\bbS)$.

Consider the likelihood function in Eq. (\ref{eq:like_phi}):
\[
f(j|\phi)=\sum_{i=1}^{W}b_{ji}\phi_{i}.\]

From the sum-to-one contraint on the parameters, it follows that $\phi_{1}=1-\sum_{i=2}^{W}\phi_{i}$.
Thus we can rewrite the previous eq. as

\begin{equation}
f(j|\phi)=b_{j1}+\sum_{i=2}^{W}(b_{ji}-b_{j1})\phi_{i}.\label{eq:loglik}\end{equation}
Hence,

\[
\frac{\partial}{\partial\phi_{k}}\log f(j|\phi)=\frac{b_{jk}-b_{j1}}{b_{j1}+\sum_{i=2}^{W}(b_{ji}-b_{j1})\phi_{i}}\qquad2<k<W.\]
From Theom. 5.1 \cite[Chapter 5]{Lehmann}, we prove that there exists a MLE that is
asymptotically efficient and asymptotically normal by showing that assumptions (A0)-(A2) and
(A)-(D) are satisfied.

\begin{pf}
(A0) Follows from~\eqref{eq:loglik}.

(A1) The support of $\phi_{i}$ for $2\leq i\leq W$ is $0<\phi_{i}<1$
subject to $\sum_{i=2}^{W}\phi_{i}\leq1$.

(A2) Observations are assumed to be independent.

(A3) Follows by the assumption that $0<\phi_{i}<1$ for $2\leq i\leq W$.

(A) We have

\[
\frac{\partial}{\partial\phi_{k}}f(j|\phi)=b_{jk},\quad2\leq k\leq W\]
and hence

\[
\frac{\partial^{3}}{\partial\phi_{m}\partial\phi_{l}\partial\phi_{k}}f(j|\phi)=0,\quad2\leq k,l,m\leq W.\]

(B) The expectation of the first logarithmic derivative of $f$ is

\begin{eqnarray*}
E_{\phi}\left[\frac{\partial}{\partial\phi_{k}}\log f(j|\phi)\right] & = & \sum_{j=1}^{W}\frac{b_{jk}-b_{j1}}{b_{j1}+\sum_{i=2}^{W}(b_{ji}-b_{j1})\phi_{i}}\left(b_{j1}+\sum_{i=2}^{W}(b_{ji}-b_{j1})\phi_{i}\right)\\
 & = & \sum_{j=1}^{W}b_{jk}-\sum_{j=1}^{W}b_{j1}\\
 & = & 1-b_{11}\\
 & = & 0.\end{eqnarray*}
As for the second derivative, we have

\begin{eqnarray*}
E\left[\frac{\partial}{\partial\phi_{l}}\log f(j|\phi)\frac{\partial}{\partial\phi_{k}}\log f(j|\phi)\right] & = & \sum_{j=1}^{W}\frac{(b_{jl}-b_{j1})(b_{jk}-b_{j1})}{\left(b_{j1}+\sum_{i=2}^{W}(b_{ji}-b_{j1})\phi_{i}\right)^{2}}\left(b_{j1}+\sum_{i=2}^{W}(b_{ji}-b_{j1})\phi_{i}\right)\\
 & = & \sum_{j=1}^{W}\frac{(b_{jl}-b_{j1})(b_{jk}-b_{j1})}{b_{j1}+\sum_{i=2}^{W}(b_{ji}-b_{j1})\phi_{i}},\end{eqnarray*}
which is equivalent to

\begin{eqnarray*}
E\left[-\frac{\partial^{2}}{\partial\phi_{l}\partial\phi_{k}}\log f(j|\phi)\right] & = & \sum_{j=1}^{W}-\left(-\frac{(b_{jk}-b_{j1})(b_{jl}-b_{j1})}{\left(b_{j1}+\sum_{i=2}^{W}(b_{ji}-b_{j1})\phi_{i}\right)^{2}}\left(b_{j1}+\sum_{i=2}^{W}(b_{ji}-b_{j1})\phi_{i}\right)\right)\\
 & = & \sum_{j=1}^{W}\frac{(b_{jl}-b_{j1})(b_{jk}-b_{j1})}{b_{j1}+\sum_{i=2}^{W}(b_{ji}-b_{j1})\phi_{i}}.\end{eqnarray*}

(C) The vectors $\left[\frac{\partial}{\partial\phi_{2}}\log f(j|\phi),\frac{\partial}{\partial\phi_{3}}\log f(j|\phi),\dots,\frac{\partial}{\partial\phi_{W}}\log f(j|\phi)\right]$
for $1<j<W$ must be linearly independent with probability 1. Note
that and $b_{jk}>0\iff j\leq k$ (in particular, $b_{j1}>0\iff j=1$).
It follows that for $j>k\geq2$

\begin{eqnarray*}
\frac{\partial}{\partial\phi_{k}}\log f(j|\phi) & = & \frac{b_{jk}-b_{j1}}{b_{j1}+\sum_{i=2}^{W}(b_{ji}-b_{j1})\phi_{i}}\\
 & = & 0,\end{eqnarray*}
whereas for $j\leq k$,

\[
\frac{\partial}{\partial\phi_{k}}\log f(j|\phi)=\frac{b_{jk}}{\sum_{i=2}^{W}(b_{ji}-b_{j1})\phi_{i}}>0.\]

Therefore, the $j-1$ leftmost entries in the $j$-th vector are 0
while the remainder are positive. Hence the vectors are linearly independent.

(D) Consider a constant $\epsilon_{j}>0$ such that $f(j|\phi)=b_{j1}+\sum_{i=2}^{W}(b_{ji}-b_{j1})\phi_{i}\geq\epsilon_{j}$
for $1\leq j\leq W$. Thus,
\begin{eqnarray*}
\left|\frac{\partial^{3}}{\partial\phi_{m}\partial\phi_{l}\partial\phi_{k}}f(j|\phi)\right| & = & \left|\frac{-(b_{jk}-b_{j1})(b_{jl}-b_{j1})\times2(b_{jm}-b_{j1})\phi_{m}(b_{j1}+\sum_{i=2}^{W}(b_{ji}-b_{j1})\phi_{i})}{\left(b_{j1}+\sum_{i=2}^{W}(b_{ji}-b_{j1})\phi_{i}\right)^{4}}\right|\\
 & = & \left|\frac{2(b_{jk}-b_{j1})(b_{jl}-b_{j1})(b_{jm}-b_{j1})\phi_{m}}{\left(b_{j1}+\sum_{i=2}^{W}(b_{ji}-b_{j1})\phi_{i}\right)^{3}}\right|\\
 & \leq & \left|\frac{2(b_{jk}-b_{j1})(b_{jl}-b_{j1})(b_{jm}-b_{j1})\phi_{m}}{\epsilon_{j}^{3}}\right|.\end{eqnarray*}
Since $M_{klm}(j)=\left|\frac{\partial^{3}}{\partial\phi_{m}\partial\phi_{l}\partial\phi_{k}}f(j|\phi)\right|<\infty$,
then $E_{\phi}[M_{klm}(j)]<\infty$ for all $k,l,m$.
\end{pf}

\section{Average set size proofs}

\begin{lemma} \label{lemma:msemphi}
Let $p$ be the sampling probability and $\hat m_\phi$ denote an unbiased
estimate of the average size of the observed sets $m_\phi$.
Then,
\[
\textrm{MSE}(\hat m_\phi) = O\left(\frac{m_{\phi}^{(2)}-m_{\phi}^{2}}{N}\right) .\]
\end{lemma}
\begin{pf}
The estimation error lower bound of the average set size is~\cite[pg.83, Proposition 3]{vanTrees}
\begin{equation}\label{eq:msemphi}
\textrm{MSE}(\hat m_\phi) \geq \frac{(1,\dots,W) (J^{(\phi)})^{-1}(1,\dots,W)^\texttt{T} - m_{\phi}^2}{N}.
\end{equation}
Lemma~\ref{lemma:inversematrix} yields
\begin{eqnarray}
&&(1,\dots,W) (J^{(\phi)})^{-1}(1,\dots,W)^\texttt{T} \nonumber \\
&=&\sum_{k=1}^{W} \sum_{i=1}^k \sum_{j=1}^k  i j  \binom{k}{j} \binom{k}{i}
\left(\frac{q}{p}\right)^{2k}
(-1)^{2k-i-j}(q^{-i}-1)(q^{-j}-1) d_k(\bphi) \nonumber \\
&=&\sum_{k=1}^{W} (q/p)^{2k}  d_k(\bphi) \left(\sum_{i=1}^k i \binom{k}{i}
        \frac{q^{-i}-1}{(-1)^{i}} \right) \left(\sum_{j=1}^k j \binom{k}{j}
        \frac{q^{-j}-1}{(-1)^{j}} \right) \nonumber \\
&=&d_1(\bphi)+\sum_{k=2}^{W} (q/p)^{2k}  d_k(\bphi) \left( \left( - \frac{1-q}{q} \right)^k \frac{k}{1-q} \right)^2 \nonumber \\
&=&\left(1-\frac{1}{p^2}\right)d_1(\bphi)+\frac{1}{p^2} \sum_{k=1}^{W}  d_k(\bphi) k^2 . \label{eq:msemphipart1}
\end{eqnarray}
Now~\eqref{eq:dj} yields
\begin{equation} \label{eq:d1}
d_1(\bphi)=\sum_{i=1}^{W} \frac{i p q^{i-1}}{1 - q^i}\phi_i
\end{equation}
and
\begin{equation*}
\begin{split}
\sum_{k=1}^{W}  d_k(\bphi) k^2&=\sum_{k=1}^{W}  \sum_{i=k}^{W} \frac{\binom{i}{k} p^k q^{i-k}}{1 - q^i}\phi_i k^2\\
&=\sum_{i=1}^{W}  \sum_{k=1}^{i} \frac{\binom{i}{k} p^k q^{i-k}}{1 - q^i}\phi_i k^2\\
&=\sum_{i=1}^{W}  \left(\sum_{k=1}^{i} \binom{i}{k} p^k q^{i-k} k^2\right) \frac{\phi_i}{1 - q^i}.\\
\end{split}
\end{equation*}
Using the relation
\[ \sum_{k=1}^{i} \binom{i}{k} x^k y^{i-k} k^2 = \left\{ \barr{ll}
  x, & i=1, \\
  ix(ix+y)(x+y)^{i-2}, & i\ge 2.
  \earr \right.
\]
yields
\begin{equation}\label{eq:msemphipart2}
\sum_{k=1}^{W}  d_k(\bphi) k^2=\sum_{i=1}^{W}  \frac{ip(ip+q)\phi_i}{1 - q^i}.
\end{equation}
Putting together~\eqref{eq:msemphi}, \eqref{eq:msemphipart1}, and \eqref{eq:msemphipart2} yields
\begin{equation}\label{eq:MSEm}
\MSE(\hat m_\phi ) \geq \left( \sum_{i=1}^W \frac{i ( p i + q^{i+1} - 2q^i +q)\phi_i }{p(1-q^i)}  - m_\phi^2  \right)/N
\end{equation}
which concludes the proof.
\end{pf}

\begin{lemma}\label{lemma:efficientavg}
Using the observed set sizes $\bbS = \{\cS_k\}_{k=1}^N$ the following
\begin{equation}\label{eq:avgestimator}
\hat{m}_\phi=\frac{\sum_{k=1}^N \cS_k}{Np}+\left(1-\frac{1}{p}\right)
    \frac{\sum_{k=1}^N {\bf 1}_{\cS_k=1}}{N},
\end{equation}
is an efficient (smallest variance) unbiased estimator of $m_\phi$.
\end{lemma}
\begin{pf}
We start by noting that
\begin{equation} \label{eq:mphi}
m_\phi=[1,...,W]\phi=[1,...,W]B^{-1}d(\bphi).
\end{equation}
Denote $z=[z_1,\ldots,z_W]=[1,...,W]B^{-1}$. From Lemma~\ref{lemma:Bstar}, we have
\begin{eqnarray}
z_i&=&\sum_{j=1}^W j b_{ji}^\star \nonumber\\
&=&\sum_{j=1}^i j \binom{i}{j} p^{-i}(-q)^{i-j}(1-q^j) \nonumber\\
&=&(-q/p)^i \sum_{j=1}^i j \binom{i}{j} \frac{1-q^j}{(-q)^j} \label{eq:zi}
\end{eqnarray}
For $i=1$ \eqref{eq:zi} yields $z_1=1$ and for $2\le i \le W$, 
\[
z_i=(-q/p)^i \left( - \frac{1-q}{q} \right)^i \frac{i}{1-q}=\frac{i}{p}.
\]
Therefore,
\[
z=\frac{[p,2, 3, \ldots, W]}{p}.
\]
Thus applying the above back into \eqref{eq:mphi} yields
\begin{equation}\label{eq:avgequation}
m_\phi=\frac{m_d}{p}+\left(1-\frac{1}{p}\right)d_1(\bphi),
\end{equation}
where $m_d=\sum_{i=1}^W i d_i$ is the expectation of average set size of observed subsets. 
Rewriting \eqref{eq:avgequation} using the set sizes $\bbS$ we get
\[
\hat{m}_\phi=\frac{1}{N} \sum_{k=1}^N
\left(\frac{\cS_k}{p}+\left(1-\frac{1}{p}\right) {\bf 1}_{\cS_k=1} \right).
\]
Based on our assumption that $\{S_k\}_{k=1}^m$ is an i.i.d.\  sequence, we have that
$\{\cS_k\}_{k=1}^N$ is also i.i.d.\ with distribution $d(\bphi)$. 
Therefore,
\[
E[\hat{m}_\phi]=E\left[\frac{\cS_k}{p}+\left(1-\frac{1}{p}\right) {\bf 1}_{\cS_k=1}\right],
\]
and
\[
\textrm{Var}[(\hat{m}_\phi)^2]=\frac{1}{N}\textrm{Var}\left[\left(\frac{\cS_k}{p}+\left(1-\frac{1}{p}\right)
        {\bf 1}_{\cS_k=1}\right)^2\right].
\]
Since
\[
E[\cS_k]=m_d=\sum_{i=1}^W i d_i(\bphi),
\]
and
\[
E[{\bf 1}_{\cS_k=1}]=d_1(\bphi),
\]
we have $E[\hat{m}_\phi]=m_\phi$ from~\eqref{eq:avgequation}, which indicates that $\hat{m}_\phi$ is unbiased.
Then
\[
E[(\cS_k)^2]=\sum_{i=1}^W i^2 d_i(\bphi),
\]
\[
E[({\bf 1}_{\cS_k=1})^2]=d_1(\bphi),
\]
and
\[
E[\cS_k {\bf 1}_{\cS_k=1}]=d_1(\bphi),
\]
yield
\begin{equation*}
\textrm{Var}[(\hat{m}_\phi)^2]=\frac{\left(1-\frac{1}{p^2}\right)d_1(\bphi)+\frac{1}{p^2} \sum_{k=1}^{W}  d_k(\bphi) k^2 - m_\phi^2}{N}.
\end{equation*}
From~\eqref{eq:msemphi} and~\eqref{eq:msemphipart1} we find that $\hat{m}_\phi$ is an unbiased estimator that achieves the Cram\'er-Rao lower bound (i.e., it is an efficient estimator).
\end{pf}

\begin{lemma} \label{lemma:msemtheta}
Let $\hat m$ denote an unbiased estimate of the average set size
$m_\theta$. Then,
\begin{eqnarray}
\textrm{MSE}(\hat m_\theta) & \geq &
\frac{1}{\eta^{2}}\Bigg(\sum_{i=1}^{W}\sum_{j=1}^{W}\frac{ij[(J^{(\bphi)})^{-1}]_{ji}}{(1-q^{j})(1-q^{i})}+
m_{\theta}^{2}\sum_{i=1}^{W}\sum_{j=1}^{W}\frac{[(J^{(\bphi)})^{-1}]_{ji}}{(1-q^{j})(1-q^{i})}-
\nonumber\\
 &  &
 2m_{\theta}\sum_{i=1}^{W}\sum_{j=1}^{W}\frac{j[(J^{(\bphi)})^{-1}]_{ji}}{(1-q^{i})(1-q^{j})}\Bigg)
    \label{eq:msemtheta}.
 \end{eqnarray}
\end{lemma}
\begin{pf}

\begin{eqnarray}
\textrm{MSE}(\hat m_\theta) & \geq & \frac{\nabla M}{\nabla\theta}\left(\frac{\nabla
        H}{\nabla\phi}(J^{(\bphi)})^{-1}\frac{\nabla H}{\nabla\phi}^{T}\right)\frac{\nabla M}{\nabla\theta}^{T}\nonumber \\
 & = & \left(\frac{\nabla M}{\nabla\theta}\frac{\nabla H}{\nabla\phi}\right)(J^{(\bphi)})^{-1}\left(\frac{\nabla M}{\nabla\theta}\frac{\nabla H}{\nabla\phi}\right)^{T}.\label{eq:l}\end{eqnarray}
where $\frac{\nabla M}{\nabla\theta}=(1,\ldots,W)$. Note that

\begin{eqnarray}
\left[\frac{\nabla M}{\nabla\theta}\frac{\nabla H}{\nabla\phi}\right]_{k} & = & \sum_{i=1}^{W}ih_{ik}\nonumber \\
 & = & \sum_{{i=1\atop i\neq k}}^{W}i\left(-\frac{\theta_{i}}{\eta(1-q^{k})}\right)+k\left(\frac{1-\theta_{k}}{\eta(1-q^{k})}\right)\nonumber \\
 & = & \frac{1}{\eta(1-q^{k})}\left(k-\sum_{i=1}^{W}i\theta_{i}\right)\nonumber \\
 & = & \frac{k-m_{\theta}}{\eta(1-q^{k})}.\label{eq:mh}\end{eqnarray}

Substituting eq. (\ref{eq:mh}) in eq. (\ref{eq:l}), we have

\begin{eqnarray*}
\textrm{MSE}(\hat m_\theta) & \geq &
\sum_{i=1}^{W}\sum_{j=1}^{W}\left(\frac{j-m_{\theta}}{\eta(1-q^{j})}\right)[(J^{(\bphi)})^{-1}]_{ji}\left(\frac{i-m_{\theta}}{\eta(1-q^{i})}\right)\\
 & = & \frac{1}{\eta^{2}}\Bigg(\sum_{i=1}^{W}\sum_{j=1}^{W}\frac{ij[(J^{(\bphi)})^{-1}]_{ji}}{(1-q^{j})(1-q^{i})}+
 m_{\theta}^{2}\sum_{i=1}^{W}\sum_{j=1}^{W}\frac{[(J^{(\bphi)})^{-1}]_{ji}}{(1-q^{j})(1-q^{i})}-\\
 &  &
 2m_{\theta}\sum_{i=1}^{W}\sum_{j=1}^{W}\frac{j[(J^{(\bphi)})^{-1}]_{ji}}{(1-q^{i})(1-q^{j})}\Bigg).\end{eqnarray*}
\end{pf}

Similarly to what we did for eq.~\eqref{eq:g_jj}, we split
eq.~\eqref{eq:msemtheta} into
three pieces to analyze its behavior.

\begin{eqnarray*}
\textrm{MSE}(\hat m_\theta) & \geq &
\frac{1}{\eta^{2}}\Bigg(\underbrace{\sum_{i=1}^{W}\sum_{j=1}^{W}\frac{ij[(J^{(\bphi)})^{-1}]_{ji}}{(1-q^{j})(1-q^{i})}}_{U_{1}}+
\underbrace{m_{\theta}^{2}\sum_{i=1}^{W}\sum_{j=1}^{W}\frac{[(J^{(\bphi)})^{-1}]_{ji}}{(1-q^{j})(1-q^{i})}}_{U_{2}}-\\
 &  &
 \underbrace{2m_{\theta}\sum_{i=1}^{W}\sum_{j=1}^{W}\frac{j[(J^{(\bphi)})^{-1}]_{ji}}{(1-q^{i})(1-q^{j})}}_{U_{3}}\Bigg).
 \end{eqnarray*}

\subsection{Analysis of $U_{1}$}

\begin{eqnarray*}
\sum_{i=1}^{W}\sum_{j=1}^{W}\frac{ij[(J^{(\bphi)})^{-1}]_{ji}}{(1-q^{j})(1-q^{i})}
 & = &
 \sum_{i=1}^{W}\sum_{j=1}^{W}\sum_{k=1}^{W}ij\binom{k}{i}\binom{k}{j}\left(\frac{q}{p}\right)^{2k}(-q)^{-i-j}d_{k}(\bphi)\\
 & = & \sum_{k=1}^{W}\left(\frac{q}{p}\right)^{2k}d_{k}(\bphi)\left(\sum_{i=1}^{k}i\binom{k}{i}\left(-q\right)^{-i}\right)^{2}\\
 & = & \sum_{k=1}^{W}\left(\frac{q}{p}\right)^{2k}d_{k}(\bphi)\left(\left(-\frac{q}{p}\right)^{-k}\frac{k}{p}\right)^{2}\qquad\textrm{using (\ref{eq:id0})}\\
 & = & \frac{1}{p^{2}}\sum_{k=1}^{W}k^{2}d_{k}(\bphi)\\
 & = & \frac{\eta}{p^{2}}\sum_{i=1}^{W}ip(ip+q)\theta_{i}\\
 & = & \eta(\sum_{i=1}^{W}i^{2}\theta_{i}+\frac{q}{p}m_{\theta}).\end{eqnarray*}

Note that $U_{1}$ is bounded by the second moment of the distribution
$\theta$.

\subsection{Analysis of $U_{2}$}

Note that $U_{2}=\frac{m_{\theta}^{2}}{\theta_{i}^{2}}A_{2}(i)$. Therefore,
we conclude that $U_{2}$ diverges if either $\theta_{W}$ decreases
exponentially in $W$ and $p<a/(a+1)$ or $\theta_{W}$ decreases
slower than exponentially in $W$ and $p<1/2$.

\subsection{Analysis of $U_{3}$}

\begin{eqnarray*}
\sum_{i=1}^{W}\sum_{j=1}^{W}\frac{j[(J^{(\bphi)})^{-1}]_{ji}}{(1-q^{i})(1-q^{j})}
 & = & \sum_{k=1}^{W}\left(\frac{q}{p}\right)^{2k}d_{k}(\bphi)\sum_{i=1}^{k}\binom{k}{i}(-q)^{-i}\sum_{j=1}^{k}j\binom{k}{j}(-q)^{-j}\\
 & = & \sum_{k=1}^{W}\left(\frac{q}{p}\right)^{2k}d_{k}(\bphi)\left(\left(-\frac{p}{q}\right)^{k}-1\right)\left(\left(-\frac{p}{q}\right)^{k}\frac{k}{p}\right)\qquad\textrm{using (\ref{eq:id1},\ref{eq:id0})}\\
 & = & \frac{1}{p}\underbrace{\sum_{k=1}^{W}kd_{k}(\bphi)}_{\eta p m_\theta}-\frac{1}{p}\underbrace{\sum_{k=1}^{W}\left(-\frac{q}{p}\right)^{k}kd_{k}(\bphi)}_{-\eta q\theta_{1}}\\
 & = & \eta(m_{\theta}+\frac{q}{p}\theta_{1}).\end{eqnarray*}

Thus,

\[
U_{3}=2m_{\theta}\eta(m_{\theta}+\frac{q}{p}\theta_{1}).\]

It is interesting to note that, counterintuitively, $U_{2}$ goes to infinity for
certain values of $p$ and $\btheta$ while $U_{1}$ and $U_{3}$ are always finite,
even though the factor $[(J^{(\phi)})^{-1}]_{ji}$ that appears inside the double summation in
$U_{2}$ is the same factor that appears multiplied by $j$ and $ji$ in $U_{1}$ and $U_{3}$,
respectively.

\subsection{Proof of Theorem~\ref{thm:average}}

Note that $U_{1}$, $U_{2}$ and $U_{3}$ are positive quantities and, moreover,
$\textrm{MSE}(\hat m_\theta)>0\Rightarrow U_{1}+U_{2}>U_{3}$. We observe that $U_{1}$ diverges
if the second moment of $\theta$ is infinite, $U_{2}$ diverges if
$\sum_{j=1}^{W}\left(\frac{q}{p}\right)^{j}\theta_{j}\rightarrow\infty$ as
$W\rightarrow\infty$, while $U_{3}$ is always finite.

\begin{pf}
\textbf{ 1) When $\theta_{W}$ decreases faster than exponentially in
$W$.}

In this case, the second moment of $\btheta$ is finite and the sum
$\sum_{j=1}^{W}\left(\frac{q}{p}\right)^{j}\theta_{j} = O(1)$ for $0 < p <
1$. Therefore, $\textrm{MSE}(m(\bbS))=O(1)$ for $0 < p < 1$.

\textbf{ 2) When $\theta_{W}$ decreases exponentially in
$W$.}

The second moment of $\btheta$ is still finite. However, we can show that the sum
$\sum_{j=1}^{W}\left(\frac{q}{p}\right)^{j}\theta_{j}$ is $\Omega(W)$ for $p \leq
a/(a+1)$ and $O(1)$ for $p > a/(a+1)$ by using an argument similar to the one used in
Section E of Appendix A. Hence, $\textrm{MSE}(m(\bbS))=\Omega(W)$ for $p \leq
a/(a+1)$ and $\textrm{MSE}(m(\bbS))=O(1)$ for $p > a/(a+1)$.

\textbf{ 3) When $\theta_{W}$ decreases more slowly than exponentially in
$W$.}


We can show that the sum $\sum_{j=1}^{W}\left(\frac{q}{p}\right)^{j}\theta_{j}$
is $\Omega(W)$ for $p < 1/2$ and $O(1)$ for $p \geq 1/2$
by using an argument similar to the one used in Section E of Appendix
A. However, the second moment of $\btheta$ shows up in $U_1$ and it can be either finite or infinite.
Although it does not affect the bound for $p <
1/2$, in which case we have $\log \textrm{MSE}(m(\bbS))=\Omega(W)$, it does
change the bound for $p \geq 1/2$. In particular, if $p = 1/2$ and $\sum_{j=1}^W j^2
\theta_j = \omega(1)$, then $\textrm{MSE}(m(\bbS))=\omega(1)$. On the other
hand, if $p = 1/2$ and $\sum_{j=1}^W j^2 \theta_j \geq O(1)$, then
$\textrm{MSE}(m(\bbS))=\Omega(1)$. Finally, if $p > 1/2$, then
$\textrm{MSE}(m(\bbS))=\Omega(1)$ as well.

\end{pf}

\section{Useful identities}

\begin{eqnarray}
\sum_{j=1}^{k}j\binom{k}{j}(-q)^{-j} & = & \left(-\frac{q}{p}\right)^{-k}\frac{k}{p}\label{eq:id0}\\
\sum_{j=1}^{k}\binom{k}{j}(-q)^{-j} & = & \left(-\frac{q}{p}\right)^{-k}-1\label{eq:id1}\\
\sum_{k=1}^{j}\binom{j}{k}\left(\frac{q}{p}\right)^{k} & = & \left(\frac{1}{p}\right)^{j}-1\label{eq:id2}\\
\sum_{k=1}^{j}\binom{j}{k}(-1)^{k} & = & -1\label{eq:id3}\\
\sum_{k=0}^{j}\binom{j}{k}(-1)^{k} & = & \begin{cases}
1 & \textrm{if \ensuremath{j=0}}\\
0 & \textrm{otherwise}\end{cases}\label{eq:id4}\end{eqnarray}

\balance
}{}

\end{document}